\documentclass[journal]{IEEEtran}
\usepackage{cite}
\usepackage{amsmath,amssymb,amsfonts}
\usepackage{graphicx}
\usepackage{textcomp}
\def\BibTeX{{\rm B\kern-.05em{\sc i\kern-.025em b}\kern-.08em
    T\kern-.1667em\lower.7ex\hbox{E}\kern-.125emX}}
\usepackage[utf8]{inputenc}
\usepackage{epsfig} 
\usepackage{mathptmx} 
\usepackage{times} 
\usepackage{bm}
\usepackage{float}
\usepackage{xcolor}

\DeclareMathOperator{\R}{\mathbb{R}}
\DeclareMathOperator{\C}{\mathbb{C}}
\DeclareMathOperator{\Pp}{\mathbb{P}}

\newcommand{\hats}[1]{\skew{-4}\hat{#1}}
\newtheorem{thm}{Theorem}[section]
\newtheorem{lem}[thm]{Lemma}

\newtheorem{cor}[thm]{Corollary}

\newtheorem{defn}{Definition}[section]

\newtheorem{rem}[thm]{Remark}
\newtheorem{assum}{Assumption}[section]

\usepackage{array}
\newcolumntype{C}[1]{>{\centering\arraybackslash}m{#1}}

\usepackage{tikz}
\usetikzlibrary{shapes,arrows}
\tikzstyle{block} = [draw, fill=white!20, rectangle, 
    minimum height=2em, minimum width=3em]
\tikzstyle{sum} = [draw, fill=white!20, circle, node distance=0.5cm]
\tikzstyle{input} = [coordinate]
\tikzstyle{output} = [coordinate]
\tikzstyle{pinstyle} = [pin edge={to-,thin,black}]
\usetikzlibrary{arrows}
\usetikzlibrary{arrows,decorations.markings}

\let\emptyset\varnothing

\usepackage{algorithm}
\usepackage{algpseudocode}
\floatname{algorithm}{Procedure}

\begin{document}

\title{On the Stability Margin and Input Delay Margin of Linear Multi-agent systems}
\author{Rajnish~Bhusal and~Kamesh~Subbarao
\thanks{The authors are with the Department of Mechanical and Aerospace Engineering, The University of Texas at Arlington, Arlington, TX 76019, USA (e-mail: rajnish.bhusal@mavs.uta.edu; subbarao@uta.edu)}
\thanks{This work was supported by the Office of Naval Research via award number N00014-18-1-2215.}}

\maketitle

\begin{abstract}

This paper provides a framework to characterize the gain margin, phase margin, and input delay margin of a linear time-invariant multi-agent system where the interaction topology is described by a graph with a directed spanning tree. The stability analysis of the multi-agent system based on the generalized Nyquist theorem is converted to finding a minimum gain positive definite Hermitian perturbation and minimum phase unitary perturbation in the feedback path of the loop transfer function. Specifically, two constrained minimization problems are solved to calculate the gain, phase and input delay margins of the multi-agent system. We further state necessary and sufficient conditions concerning stability of the multi-agent system independent of gain and phase perturbations, and input delay.
\end{abstract}



\section{Introduction}
\label{sec:introduction}
In recent years, a significant amount of research efforts have been focused on the cooperative control of multi-agent systems and a variety of distributed control protocols have been proposed to perform desired cooperative tasks among the agents. The distributed control protocols for multi-agent systems find their applications in formation control \cite{dong2016time}, flocking \cite{jafari2019biologically}, rendezvous of unmanned aerial vehicles \cite{zhang2015cooperative}, attitude synchronization among multiple spacecrafts \cite{cai2014leader}, among others. 

With increasing applications, stability and robustness-based analysis of multi-agent systems has also drawn significant attention. This paper in particular, provides a framework for calculation of stability margin and input delay margin for a group of multiple agents in the networked interconnection. For a single-input single-output (SISO) system, classical input-output stability criteria based on Nyquist, Popov and circle theorems aid to characterize the allowable gain and phase variation (stability margin) in the loop at each frequency and tolerable limits of open-loop modeling errors. Generalizations of the aforementioned theorems to multi-input multi-output (MIMO) systems is not straightforward, and several works such as \cite{safonov1981multiloop, lehtomaki1981robustness, nie2010exact, wang2008loop} suitably characterize the MIMO stability margins. In the context of multi-agent systems, the stability margin serves as a robustness measure against gain and phase variations for the group of agents. A networked multi-agent system is a multiloop feedback system and with suitable analysis, the aforementioned works to characterize the multiloop stability margin can be extended to the context of multi-agent systems. On that note, Tonetti and Murray \cite{tonetti2010limits} have considered disturbance rejection based graph topology-design strategies for multi-agent systems by calculating the gain and phase margins of interconnected systems upon analyzing the networked sensitivity function matrix. However, the analysis in \cite{tonetti2010limits} stems from the assumption that each individual agent is a SISO system. In  \cite{gattami2004frequency}, although a Nyquist-like criterion is presented to analyze stability of the interconnected system of agents, the stability margins of the interconnected system are not characterized. Also, Kim \cite{kim2017stability} characterized the stability margin of SISO multi-agent systems based on the minimum singular value of the loop transfer function matrix.

Moreover, multi-agent systems need to exchange information among agents over a communication network, which invariably, is prone to time delays. The presence of time delay may significantly degrade closed-loop performance, and even cause instability. As mentioned in \cite{cao2013overview}, two types of time delays, input delay and communication delay, have been considered in the literature. Input delay is related to processing and connecting time for the packets arriving at each agent while communication delay refers to the time for transferring information between agents. As discussed in \cite{tian2008consensus}, for integrator dynamics, when certain connectivity condition is satisfied by the topology graph, the
consensusability conditions are independent of communication delays, but dependent on input delays. Therefore stability criteria for multi-agent systems with input delays have been attracting great attention over the years \cite{xu2013input, zhang2018synchronization}. For integrator dynamics of agents, the time delay problem has been discussed in \cite{olfati2004consensus}, which provides necessary and sufficient conditions for the maximum delay such that the multi-agent system reaches consensus from arbitrary initial conditions. Stability conditions in terms of linear matrix inequalities (LMIs) using Lyapunov Krasovskii techniques for single integrator dynamics of agents under consensus protocol with input delays are provided in \cite{lin2008average}. In \cite{munz2010delay}, robust consensus conditions for multi-agent system consisting of SISO agents in undirected network subject to heterogeneous feedback delays are derived from frequency-dependent and delay-dependent convex sets. Furthermore in \cite{xu2013input}, the input delay margin for consensus among agents under undirected graph topology with scalar dynamics and single input vector dynamics with a single unstable open-loop pole is derived. Recently in \cite{zhang2019state}, static consensus protocols under undirected graph topology have been derived for multi-agent systems with nonuniform input delays. Although, most of the works in the literature for high-order multi-agent systems with input delay are restricted to undirected graphs, some of the recent works for multi-agent systems with input delay under directed graph topology can be found in \cite{zhang2017distributed, zhao2017guaranteed}.

In this paper, the problems of destabilizing gain, phase, and input delays applied to multi-agent systems consider a group of agents modeled as high-order linear dynamical systems. The interconnections within the group are modeled using a graph with at least one directed spanning tree. We develop a framework to characterize the stability margins as a direct multivariable generalization of the complex units used in SISO phase analysis. More specifically, we are concerned with the stability of the collective dynamics of the agents subjected to complex perturbations. The application of such perturbation analysis is significant in the areas where any errors such as signal interference or time delays in sensors introduce significant gain and phase shifts which might affect the collective stability of networked agents. The overall effect of such errors can be modeled as a complex perturbation in the feedback loop \cite{kim2017stability}. On the other hand, it is well-known that frequency based representation of the time delay links it with the phase lag in the system with no gain change. This motivates us to obtain the input delay margin of a multi-agent system based on the unitary phase perturbation of the system's loop transfer function in the feedback path. In this paper, we consider the delays in the input of all the agents to be uniform. The work carried out in this paper uses some of the results from the work carried on multivariable gain and phase margins in \cite{bar1990phase} and \cite{bar1991multivariable}. The major contributions of the paper can be enumerated as follows:
\begin{itemize}
\item[(i)] With a controller that guarantees the closed loop stability of a high order (linear) system, we transform the stability criteria for consensus among $N$ identical agents with a distributed control protocol (using generalized Nyquist's criteria) to an equivalent stability criteria of $N-1$ MIMO loop transfer functions. 

\item[(ii)] We develop a unified framework to compute the gain margin, phase margin and input delay margin for multi-agent systems to achieve consensus. The problem of calculating the stability margins and input delay margin is converted to finding eigenvalues of multiplicative perturbation in the feedback paths of a set of multi-input multi-output (MIMO) loop transfer functions which involves solving a constrained minimization problem. We do not impose any restrictions on the dynamics of agents and on the graph topology, except that the graph structure should have atleast a directed spanning tree which is imperative for consensus.

\item[(ii)] The closed loop stability of a general MIMO system independent of gain and phase perturbations, and input delay can be treated as a robust stability problem and suitable small gain conditions can be derived for the stability. To that end, we develop necessary and sufficient conditions for gain-independent, phase-independent and delay-independent stability of multi-agent systems which can be considered to be extended small gain conditions.
\end{itemize}    

The paper is organized as follows. In section \ref{sec:primer}, we briefly review principal concepts of graph theory and formulate the problem in consideration. Section \ref{sec:stability_general} discusses the stability of multi-agent system with or without input delay in general. Main results of the paper are presented in section \ref{gm_pm_tdm}. Numerical examples are presented in section \ref{sec:results} and the conclusions of the paper are reported in section \ref{sec:conclusion}.


\section{Preliminaries and Problem Formulation}
\label{sec:primer}
\subsection{Notations}
For a vector $\mathbf{x} \in \R^n$, $\|\mathbf{x}\|$ denotes its 2-norm. For a $\xi \in \C$, its real part is denoted by $\text{Re}(\xi)$ and its imaginary part by $\text{Im}(\xi)$.  For a matrix $\mathbf{T} \in \R^{n \times n}$, $\lambda_i(\mathbf{T})$, $i=1,2, \dots, n$, denote its eigenvalues, $\sigma_i(\mathbf{T})$, $i=1,2, \dots, n$, denote its singular values and $\text{det}(\mathbf{T})$ denotes its determinant. We denote complex conjugate transpose of a complex matrix $\mathbf{T} \in \C^{n \times n}$ by $\mathbf{T}^*$. $\mathbf{A} \otimes \mathbf{B}$ denotes the Kronecker product of matrices $\mathbf{A}$ and $\mathbf{B}$. $\mathbf{1}_n$ denotes a $n$-dimensional vector of ones; $\mathbf{I}_n$ denotes the identity matrix of dimension $n \times n$. For two real (complex) vectors $\mathbf{x}$ and $\mathbf{y}$, $\left\langle\mathbf{x},\mathbf{y}\right\rangle$ denotes their inner product (Hermitian inner product). We denote $\text{col}(\mathbf{x}_1, \mathbf{x}_2, \dots, \mathbf{x}_n)$ as concatenation of vectors $\mathbf{x}_1, \mathbf{x}_2, \dots, \mathbf{x}_n$ such that $\text{col}(\mathbf{x}_1, \mathbf{x}_2, \dots, \mathbf{x}_n) = [\mathbf{x}^\text{T}_1, \mathbf{x}^\text{T}_2, \dots, \mathbf{x}^\text{T}_n ]^\text{T}$.

\subsection{Algebraic Graph Theory}
Denote $\mathcal{G} = (\mathcal{V},\mathcal{E},\mathcal{A})$ as a weighted graph composed of a set of nodes $\mathcal{V} = \left\lbrace 1, 2, \dots, N \right\rbrace$ and a set of ordered pairs of nodes, called edges $\mathcal{E}=\left\lbrace (i,k)  | \ i,k \in \mathcal{V} \right\rbrace \subseteq \mathcal{V} \times \mathcal{V}$, and the adjacency matrix $\mathcal{A} = [a_{ik}] \in \R^{N \times N}$. Each of the edges of a graph $(i,k)$ is associated with a non-negative weight $a_{ik}$ such that  $a_{ik} > 0$ if $(k, i$) $\in$ $\mathcal{E}$ and $a_{ik} = 0$, otherwise. Node $k$ is the neighbor of $i$ if  $(k, i) \in \mathcal{E}$ and the set of neighbors of node $i$ can be represented as $\mathcal{N}_i$. The Laplacian matrix of a graph $\mathbf{L}=[l_{ik}] \in \R^{N \times N}$ is defined as $l_{ii} = \sum_{k \neq i} a_{ik}$ and $l_{ik} = -a_{ik}$, where $i \neq k$.

The graph $\mathcal{G}$ is said to be strongly connected if $i, k$ are connected for all distinct nodes $i, k \in \mathcal{V}$. Directed path from node $i$ to node $k$ is defined as a sequence of successive edges in the form $\{(i, l), (l, m), \dots, (n, k) \}$.  A root $r$ is a node such that for each node $i$ different from $r$, there is a directed path from $r$ to $i$. A directed tree is a directed graph, in which there is exactly one root and every node except for this root itself has exactly one parent. A directed spanning tree is a directed tree consisting of all the nodes and some edges in $\mathcal{G}$. A directed graph contains a directed spanning tree if one of its subgraphs is a directed spanning tree  \cite{yu2011second}.

\begin{assum} \label{connected_root_node}
Throughout the paper, the graph is assumed to be strongly connected with atleast one directed spanning tree.
\end{assum}

\begin{lem}\cite{ren2007information, olfati2007consensus}
Let $\mathcal{G}$ be a strongly connected graph with atleast one directed spanning tree. Let $\lambda_i$, $i=1,2, \dots, N$ be the eigenvalues of the Laplacian matrix. Then, $\lambda_1=0$ is always a simple and the smallest eigenvalue of the Laplacian matrix, and $\text{Re}(\lambda_k)>0$, for all $k=2, \dots, N$.
\end{lem}


\subsection{Problem Formulation}
\subsubsection{Multi-Agent System without delay}
\label{sec:mas_gen}
Consider a group of $N$ identical agents. The dynamics of the $i$th agent is described by the following linear time-invariant (LTI) system
\begin{equation}
\dot{\mathbf{x}}_i(t) = \mathbf{A} \mathbf{x}_i(t) + \mathbf{B} \mathbf{u}_i(t), \qquad i=1,\dots, N \label{single_agent_dyn}
\end{equation}
where $\mathbf{A} \in \R^{n \times n}$, $\mathbf{B} \in \R^{n \times m}$ are the system matrices with $\mathbf{x}_i \in \R^n$ as the state and $\mathbf{u}_i \in \R^m$ as the input of the $i$th agent. The LTI continuous dynamics of each agent can also be represented by the loop transfer function in frequency domain as
\begin{equation}
\mathbf{P}(s) = (s \mathbf{I}_n-\mathbf{A})^{-1} \mathbf{B}
\end{equation}
which is the linear mapping of Laplace transform from the input $\mathbf{u}_i(t)$ to the state $\mathbf{x}_i(t)$.

\begin{defn}
The group of agents are said to reach consensus under any control protocol $\mathbf{u}_i$ if for any set of initial conditions $\{\mathbf{x}_i(0)\}$ there exists $\mathbf{x}^c\in \R^n$ such that $\lim_{t \to \infty} \mathbf{x}_i(t)=\mathbf{x}^c$ for all $i=1, 2, \dots, N$.
\end{defn}

\begin{assum} \label{assum_stab}
$(\mathbf{A}, \mathbf{B})$ is stabilizable.
\end{assum}

With assumption \ref{assum_stab}, let each of the agents $i=1, 2, \dots, N$ have identical feedback controller $\mathbf{K} \in \R^{m \times n}$ such that $\mathbf{A}-\mathbf{B}\mathbf{K}$ is stable. For the $i$th agent with plant transfer function $\mathbf{P}(s)$ and a state feedback controller $\mathbf{K}(s)$, we define $\mathbf{H}(s) \in \R^{n \times n}$ to be the loop transfer function as seen when breaking the loop at the output of the plant. Thus, for each agent $i=1,2, \dots, N$, we have
\begin{equation}
\mathbf{H}(s) = \mathbf{P}(s) \mathbf{K}(s). \label{H_s_tf}
\end{equation}

We consider following static distributed control protocol based on the relative states between neighboring agents as discussed in \cite{li2009consensus, li2014distributed, zhang2011optimal}:
\begin{equation}
\label{control_protocol}
\mathbf{u}_i(t) =  c\mathbf{K} \sum_{k\in \mathcal{N}_i} a_{ik} (\mathbf{x}_k (t) - \mathbf{x}_i(t)), \qquad i = 1,2, \dots, N.
\end{equation}
where $c$ is the coupling gain and $\mathbf{K}$ is the feedback gain matrix. The approach to calculate $c$ would be discussed later in Section \ref{sec:stability_general}. With control protocol in \eqref{control_protocol}, the overall global closed-loop dynamics can be written as
\begin{equation}
\dot{\mathbf{x}}(t) = (\mathbf{I}_N \otimes \mathbf{A}) \mathbf{x}(t) - c(\mathbf{L} \otimes \mathbf{B} \mathbf{K}) \mathbf{x}(t) \label{closed_loop}
\end{equation}
where $\mathbf{x} = [\mathbf{x}^T_1, \dots, \mathbf{x}^T_N]^T \in R^{Nn}$ is the global state of multi-agent system. Now, the overall loop transfer function of multi-agent system is $\mathbf{G}(s) = \mathbf{\hat{H}}(s) \mathbf{\hats{L}}$, where $\mathbf{\hat{H}}(s) = \mathbf{I}_N \otimes \mathbf{H}(s)$ and $\mathbf{\hats{L}}=c \left(\mathbf{L} \otimes \mathbf{I}_n \right)$.

In this paper, we intend to characterize the gain and phase margin of the closed-loop system \eqref{closed_loop} with state feedback controller $\mathbf{K}$.

\subsubsection{Multi-agent system with input delay}
\label{sec:mas_td}
\begin{figure}[h]
\centering
\begin{tikzpicture}[auto, node distance=2.5cm,>=latex']
    \node [block, name=delay] {$e^{-s \tau}$};
    \node [block, right of=delay,
            node distance=2.5cm] (system) {$\mathbf{P}(s)$};
    \node [input, left of=delay] (input) {};
    \draw [->] (input) -- node {$\mathbf{u}_i(t)$} (delay);
    \draw [->] (delay) -- node[name=u] {$\mathbf{u}_i(t-\tau)$} (system);
    \node [output, right of=system] (output) {};
    \draw [->] (system) -- node [name=x] {$\mathbf{x}_i(t)$}(output);
        node [near end] {} (sum);
\end{tikzpicture}
\caption{Schematic representation of input delay for $i$th agent} \label{inputdelay_config}
\end{figure}
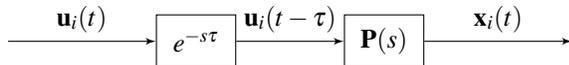

Let us now consider a problem of multi-agent system with $N$ agents subjected to input delay. We assume the input delays to be uniform for all agents. The dynamics of $i$th agent in the presence of input delay can be written as
\begin{equation}
\mathbf{x}_i (t) = \mathbf{A} \mathbf{x}_i(t) + \mathbf{B} \mathbf{u}_i (t-\tau), \quad i = 1,2, \dots, N \label{gen_dynamics_single_tdelay}
\end{equation}
where $\tau$ is the delay in the input of the agents. Figure \ref{inputdelay_config} illustrates the schematics of input delay for $i$th agent.

The presence of input delay governs the multi-agent system such that each agent $i$, for all $i = 1, 2, \dots, N$ receives the state information of its neighbor and its own state information with a delay of $\tau$.  Let Assumption \ref{assum_stab} holds for \eqref{gen_dynamics_single_tdelay}. With the distributed control protocol in \eqref{control_protocol}, the closed loop dynamics of $i$th agent can be written as
\begin{equation}
\label{closed_loop_withdelay_agent_i}
\mathbf{x}_i (t)  = \mathbf{A} \mathbf{x}_i(t) + c \mathbf{B} \mathbf{K} \left(\sum_{k\in \mathcal{N}_i} a_{ik} \left(\mathbf{x}_k (t-\tau) - \mathbf{x}_i(t-\tau)\right) \right)
\end{equation}

With $\mathbf{x} = [\mathbf{x}^T_1, \dots, \mathbf{x}^T_N]^T \in R^{Nn}$ as the global state of multi-agent system, the overall global closed-loop dynamics for input delay multi-agent system can be written as,
\begin{equation}
\dot{\mathbf{x}}(t) = (\mathbf{I}_N \otimes \mathbf{A}) \mathbf{x}(t) - c(\mathbf{L} \otimes \mathbf{B} \mathbf{K}) \mathbf{x}(t-\tau). \label{closed_loop_withdelay}
\end{equation}

For the multi-agent system with input delay, we are interested in finding the input delay margin $\tau^*$ such that the global closed-loop system \eqref{closed_loop_withdelay} is stable for any $\tau \in [0, \tau^*]$. 


\section{Stability in Multi-agent Systems for Consensus}
\label{sec:stability_general}
In this section, we discuss the stability conditions required for multi-agent systems to reach the consensus.

\subsection{Consensus in multi-agent systems without delay}

\begin{lem} 
\label{lem_alternate_stability}
 If $\mathbf{\hat{H}}(s) \mathbf{\hats{L}}$ has $p_u$ unstable poles, the closed loop system \eqref{closed_loop}  is stable, iff any of the following two statements hold:
\begin{itemize}
\item[(1)] The Nyquist plot of det[$\mathbf{I}_{nN} + \mathbf{\hat{H}}(s) \mathbf{\hats{L}}$] makes $p_u$ anti-clockwise encirclements of the origin. 
\item[(2)] The Nyquist plot of $\prod^N_{p=2} \text{det}[\mathbf{I}_{n} + c\lambda_p \mathbf{H}(s)]$ makes $p_u$ anti-clockwise encirclements of the origin; where $\{\lambda_p \}^N_{p=1}$ are the eigenvalues of $\mathbf{L}$.
\end{itemize}
\end{lem}

\begin{IEEEproof}
The statement (1) is the direct consequence of generalized Nyquist Theorem for the closed loop stability of a multi-agent system. Now, the equivalence of the above two statements can be shown with the help of Schur decomposition of Laplacian matrix $\mathbf{L}$ as $ \mathbf{L} = \mathbf{S} \mathbf{T} \mathbf{S}^{*}$, where $\mathbf{S}$ is a unitary matrix and $\mathbf{T}$ is an upper triangular matrix. Since $\mathbf{T}$ is an upper triangular matrix with same spectrum as $\mathbf{L}$, the eigenvalues of $\mathbf{L}$ are the diagonal entries of $\mathbf{T}$. Moreover, $\mathbf{T}$ can be decomposed as,
\begin{equation}
\mathbf{T} = \bm{\Lambda} + \bm{\Gamma}. \label{eq_T}
\end{equation}  
where $\bm{\Lambda}$ is a diagonal matrix consisting of eigenvalues $\{\lambda_p \}^N_{p=1}$ of $\mathbf{L}$ and $\bm{\Gamma}$ is a strictly upper triangular matrix. Since, $\mathbf{\hats{L}} = c\left(\mathbf{L} \otimes \mathbf{I}_n\right)$, one can write 
\begin{equation*}
\begin{aligned}
\text{det} \left[\mathbf{I}_{nN} + \mathbf{\hat{H}} (s) \mathbf{\hats{L}} \right] =\ &\text{det}\left[\mathbf{I}_{nN} + c\mathbf{\hat{H}}(s) \left(\mathbf{S} \otimes \mathbf{I}_n\right) \mathbf{\hat{T}} \left(\mathbf{S}^* \otimes \mathbf{I}_n \right)\right]\\
=\ & \text{det}\left[\left(\mathbf{S} \otimes \mathbf{I}_n\right)(\mathbf{I}_{nN} +  c\mathbf{\hat{H}}(s) \mathbf{\hat{T}}) \left(\mathbf{S}^* \otimes \mathbf{I}_n \right)\right]\\
=\ &\text{det}[\mathbf{I}_{nN} +  c\mathbf{\hat{H}}(s)  \mathbf{\hat{T}}].
\end{aligned}
\end{equation*}
where $\mathbf{\hat{T}} = \mathbf{T} \otimes \mathbf{I}_n$. As $\mathbf{\hat{H}}(s)$ is block diagonal and $\mathbf{\hat{T}}$ is block upper triangular, one can write
\begin{equation*}
\begin{aligned}
\text{det}[\mathbf{I}_{nN} + c\mathbf{\hat{H}}(s)  \mathbf{\hat{T}}]
=\ & \text{det}[\mathbf{I}_{nN} + c \mathbf{\hat{H}}(s)  \left(\bm{\Lambda} \otimes \mathbf{I}_n \right)]\\ 
=\ & \prod^N_{p=1} \text{det}[\mathbf{I}_{n} +  c\lambda_p \mathbf{H}(s)]\\
=\ & \prod^N_{p=2} \text{det}[\mathbf{I}_{n} +  c\lambda_p \mathbf{H}(s)]
\end{aligned}
\end{equation*}
The last equality comes from the fact that $\lambda_1 = 0$.
\end{IEEEproof}

\begin{rem}
\label{rem_alternate_stability}
Lemma \ref{lem_alternate_stability} implies that stability of multi-agent system is equivalent to the stability of following $p$ transformed systems
\begin{equation}
\dot{\bm{\xi}}_p(t) = \mathbf{A} \bm{\xi}_p(t) + \mathbf{B} \mathbf{u}_p(t), \quad \forall \ p=2, 3, \dots, N \label{perturbed_trans_leaderless}
\end{equation}
where, $\bm{\xi}_p$ is the state vector and $\mathbf{u}_p(t)$ is the input of the $p^\text{th}$ system which is given by  $\mathbf{u}_p(t)= -\bar{\mathbf{K}}_p \bm{\xi}_p(t)$ with $\bar{\mathbf{K}}_p =  c\lambda_p \mathbf{B}\mathbf{K}$. The essence of Lemma \ref{lem_alternate_stability} is similar to the discussion carried out for formation control of multi-agent systems in \cite{olfati2007consensus} where the authors conclude that if the controller $\mathbf{K}$ stabilizes the transformed system for all $\lambda_p$ other than the zero eigenvalue, it stabilizes the relative dynamics of formation. Alike in \eqref{H_s_tf}, we define the loop transfer functions of the transformed systems as
\begin{equation}
\mathbf{G}_p (s) = \mathbf{P}(s) \bar{\mathbf{K}}_p(s) = (s \mathbf{I}_n-\mathbf{A})^{-1} \mathbf{B}  c \lambda_p \mathbf{K}(s), \quad \forall p=2, 3, \dots, N. \label{nominal_loop}
\end{equation}
\end{rem}

\subsubsection{Selection of $\mathbf{K}$ and $\mathbf{c}$}
As stated earlier, $\mathbf{K}$ is selected such that the dynamics of the individual agent is stable before the interconnection, i.e., $\mathbf{A}-\mathbf{B}\mathbf{K}$ is Hurwitz. Now, the value of $c$ is selected such that the consensus among the agents is achieved, i.e. $\mathbf{A}- \mathbf{B} \bar{\mathbf{K}}_p= \mathbf{A}- c \lambda_p \mathbf{B} \mathbf{K}$ for $p=2, 3, \dots, N$ are Hurwitz, where $\lambda_p$ are the eigenvalues of Laplacian matrix. In this paper, we select $c$ based on the consensus region approach discussed in \cite{li2014distributed}. The consensus region of a multi-agent system can be defined as $\mathcal{S}(\sigma) = \{\sigma \in \C \mid \mathbf{A}- \sigma \mathbf{B} \mathbf{K} \text{ is Hurwitz} \}$. From \cite{li2014distributed}, for the agents to reach consensus, the coupling gain $c$ is to be selected such that $c \lambda_p \in \mathcal{S}(\sigma)$.

We make following assumption throughout the paper for further analysis.
\begin{assum}
\label{assum_trans_stability}
$\mathbf{A}- \mathbf{B} \bar{\mathbf{K}}_p$ is stable, for all $p = 2,3, \dots, N$
\end{assum}

\subsection{Consensus in multi-agent systems with input delay}
\label{sec:mas_td}
The loop transfer function of the multi-agent system \eqref{closed_loop_withdelay} can be written as $\mathbf{\hat{H}}(s) \mathbf{\hats{L}}  e^{-s \tau}$. Now, Lemma \ref{lem_alternate_stability} can be extended for the multi-agent system with delay and the stability of multi-agent system with input delay \eqref{closed_loop_withdelay} is equivalent to the stability of following $p$ transformed systems 
\begin{equation}
\dot{\bm{\xi}}_p(t) = \mathbf{A} \bm{\xi}_p(t) + \mathbf{B} \mathbf{u}_p(t-\tau), \quad \forall \ p=2, 3, \dots, N \label{perturbed_trans_tdelay}
\end{equation}
where $\bm{\xi}_p$ is the state vector and $\mathbf{u}_p(t-\tau) = -\bar{\mathbf{K}}_p \bm{\xi}_p(t-\tau)$ with $\bar{\mathbf{K}}_p =  c\lambda_p \mathbf{B}\mathbf{K}$,  is the delayed input of the $p^\text{th}$ system. Moreover, the loop transfer function of the transformed system \eqref{perturbed_trans_tdelay} becomes $(s \mathbf{I}_n-\mathbf{A})^{-1} \mathbf{B} c \lambda_p \mathbf{K}(s) e^{-s \tau}$. 


\section{Stability Margins and Input Delay Margin of Multi-agent System}
\label{gm_pm_tdm}
The stability margin serves as a robustness measure against gain and phase variations in the feedback path of the group of agents. Moreover, as stated earlier, time delays in multi-agent systems are practically unavoidable. In this section, we provide a computational framework to characterize the stability margins, namely gain and phase margins of the delay-free system \eqref{closed_loop}, and input delay margin of multi-agent system with input delay \eqref{closed_loop_withdelay}. 

\begin{defn} \cite{gantmakher1959theory}
The polar decomposition of a matrix $\mathbf{T} \in \C^{r \times t}$ with $r \geq t$ can be written as $\mathbf{T} = \mathbf{R} \mathbf{U}$ where $\mathbf{R} \in \C^{r \times t}$ is a positive semi-definite Hermitian matrix and $\mathbf{U} \in \C^{t \times t}$ is a unitary matrix. \label{defn_polar}
\end{defn}

In this paper, we calculate the stability margins and the input delay margin of the multi-agent system by assessing the characteristics of the perturbed loop transfer function $\mathbf{G}_p(j \omega_p) \bm{\Delta}_p$ where $\bm{\Delta}_p \in \C^{n \times n}$ is the multiplicative complex perturbation in the feedback path of the loop transfer function $\mathbf{G}_p (s)$. Here we consider different mathematical structures of $\bm{\Delta}_p$ depending upon the type of margin that is being computed, i.e., for computation of gain margin, phase margin and input delay margin, $\bm{\Delta}_p$ would be complex gain, phase and delay perturbations, respectively. 

The polar decomposition is a generalization to complex matrices of the familiar polar representation $z = r e^{j \phi}, r \geq 0$ of a complex number $z \in \C$. From Definition \ref{defn_polar} we can polar-decompose $\bm{\Delta}_p$ as, $\bm{\Delta}_p = \mathbf{R} \mathbf{U}$. The unitary factor  $e^{j \phi}$ of $z$ corresponds to unitary matrix $\mathbf{U} = e^{\bm{\Sigma}_p}$ of $\bm{\Delta}_p$, where $\bm{\Sigma}_p$ is a skew Hermitian matrix with phase information of $\mathrm{\bm{\Delta}}_p$ and $r = |z|$ of $z$ corresponds to the Hermitian factor $\mathbf{R}$ of $\bm{\Delta}_p$ \cite{zielinski1995polar}. We assume that the complex perturbation $\bm{\Delta}_p$ is nonsingular and thus, the polar decomposition is unique and $\mathbf{R}$ is positive definite Hermitian.

\begin{defn}
The complex perturbation $\bm{\Delta}_p$ for any $p=2, 3, \dots, N$ in the loop transfer function $\mathbf{G}_p(s)$ is said to be destabilizing at frequency $\omega_p \in \R$ if 
\begin{equation}
\text{det} (\mathbf{I} + \mathbf{G}_p(j \omega_p) \bm{\Delta}_p) = 0 \label{defn_destabilizing_perturbation}
\end{equation}
\end{defn}

\begin{lem}
\label{lem_instability_overall}
If there exists a destabilizing $\bm{\Delta}_p$ in the feedback path of $\mathbf{G}_p(j \omega_p)$ for any $p=2, 3, \dots, N$, the original loop transfer function $\mathbf{G}(s) = \hat{\mathbf{H}} \mathbf{\hats{L}}$ becomes unstable.
\end{lem}

\begin{IEEEproof}
From Lemma \ref{lem_alternate_stability} and Remark \ref{rem_alternate_stability} the stability of original loop transfer function is equivalent to the stability of $p$ transformed loop transfer functions $\mathbf{G}_p(s)$ simultaneously. Thus, if there exists a unitary $\bm{\Delta}_p$ that satisfies \eqref{defn_destabilizing_perturbation} for any $p=2, 3, \dots, N$, it destabilizes the $p^{th}$ transformed system and equivalently, the original loop transfer function $\mathbf{G}(s) = \hat{\mathbf{H}} \mathbf{\hats{L}}$.
\end{IEEEproof}

As stated earlier, input delay margin can be associated with stabilizing ranges of phase in the system which motivates us to compute the input delay margin by considering the phase perturbations in the system. Thus, we first provide a framework to characterize phase and input delay margins in a consecutive manner and provide a framework to compute gain margin separately.

\subsection{Phase Margin and Input Delay Margin}
\label{sec:pm_td}
In order to characterize the phase margin and input delay margin, it is assumed that $\mathbf{R}$ is lumped into the loop transfer function or assumed to be an identity matrix. Thus the analysis presented in the paper for characterizing phase margin and delay margin considers $ \mathrm{\bm{\Delta}}_p = \mathbf{U}=e^{\bm{\Sigma}_p}$. Hereafter, we use $ \mathrm{\bm{\Delta}}_p$, $\mathbf{U}$, $e^{\bm{\Sigma}_p}$ would be used interchangeably for characterizing phase and input delay margins. The following Lemma is an extension to the work carried out by Wang et al. in \cite{wang2008loop}, wherein $\bm{\Delta}_p$ was assumed to be structured diagonal perturbation; however in this work, we consider phase perturbations to be in the entire set of unitary matrices and not necessarily to be diagonal.

\begin{lem} \label{lem_stab_boundary}
The stabilizing boundary of phase is symmetric with respect to the origin. 
\end{lem}

\begin{IEEEproof}
Let us start by saying $(\phi_1, \phi_2, \dots, \phi_{n})$ is the point on the stabilizing boundary, then there exists some critical frequency $\omega_{c_p}$ for all $p=2, 3, \dots, N$ such that
\begin{equation*}
\text{det} [\mathbf{I} + \mathbf{G}_p(j \omega_{c_p}) \mathrm{\bm{\Delta}}_p] =\text{det} [\mathbf{I} + \mathbf{G}_p(j \omega_{c_p}) e^{\bm{\Sigma}_p}] = 0
\end{equation*}
As stated before, $\bm{\Sigma}_p$ is a skew Hermitian matrix with $\bm{\Sigma}_p = -\bm{\Sigma}^*_p$. The eigenvalue decomposition of $\bm{\Sigma}_p$ can be written as $\bm{\Sigma}_p = \mathbf{P} \Lambda_{\Delta} \mathbf{P}^*$, where $\mathbf{P}$ is a unitary matrix of eigenvectors and $\Lambda_{\Delta}$ is a diagonal matrix of eigenvalues of $\bm{\Sigma}_p$. As the phase information of $\mathrm{\bm{\Delta}}_p$ is contained in the unitary matrix $\mathbf{U}$ of the polar decomposition, the eigenvalues of $\mathbf{U}$ all lie on the unit circle such that, $\lambda_k(\mathbf{U}) = e^{j \phi_k}$ for all $k=1, 2, \dots, n$. Moreover, $\text{Im} \{\lambda_k(\bm{\Sigma}_p)\} = \text{arg} \{\lambda_k(\mathbf{U}) \} = \phi_k$ which implies $\Lambda_{\Delta} = \text{diag}(j \phi_1, j\phi_2, \dots, j\phi_{n})$. Clearly, one can write
\begin{equation}
\text{det} [\mathbf{I} + \mathbf{G}_p (j \omega_{c_p}) \mathbf{P} e^{ \{\text{diag}(j \phi_1, j\phi_2, \dots, j\phi_{n}) \}} \mathbf{P}^*] = 0 \label{eq_stab_bound}
\end{equation} 

On taking conjugate on the both sides of \eqref{eq_stab_bound}, we get
\begin{equation*}
\text{det} [\mathbf{I} + \mathbf{G}_p (-j \omega_{c_p}) \mathbf{P}^* e^{ \{\text{diag}(-j \phi_1, -j\phi_2, \dots, -j\phi_{n}) \}} \mathbf{P}]= 0
\end{equation*}

Thus, it can be asserted that for the point $(-\phi_1, -\phi_2, \dots, -\phi_{n})$, there exists an $-\omega_{c_p}$ such that the closed-loop system is marginally stable. This implies that $(-\phi_1, -\phi_2, \dots, -\phi_{n})$ is also the point on the stabilizing boundary. 
\end{IEEEproof}

\begin{rem}
\label{rem_stab_boundary}
By Lemma \ref{lem_stab_boundary}, the stabilizing borders of loop phases are symmetric with respect to the origin, the values of $\omega_{c_p}$ are also symmetric with respect to the origin. This property hints that one only needs to examine the frequency response for nonnegative frequencies, while the analysis for the other half of the frequency range follows that of nonnegative frequency range due to symmetry. This simplification is analogous to the analysis of half-sectorial systems in the work of Chen et al. \cite{chen2019phase}. Further, $e^{j \phi_k}$ is a periodic function in $\phi_k$ with a period of $2 \phi_k$ and thus, one only needs to consider $\phi_k \in (-\pi, \pi]$ and for discussing stability, it can be further narrowed to $\phi_k \in [0, \pi]$, for all $k=1, 2, \dots, n$. Moreover, phases of $\bm{\Delta}_p = \mathbf{U}$ for all $p = 2, \dots, N$ can be calculated as $\phi_k = |\text{Im} \{\lambda_k(\mathbf{U}) \}|$ for all $k=1, 2, \dots, n$ and in turn phase of $\bm{\Delta}_p$ can be defined as $\max(|\text{Im}(\lambda_k(\mathbf{U})|))$ in $[0, \pi]$.
\end{rem}

\subsubsection{Stability of multi-agent system independent of unitary phase perturbations}
We provide following necessary and sufficient conditions such that the multi-agent system is stable for any unitary phase perturbation in the feedback path. These conditions can be considered to be an extended small gain conditions in robust stability analysis. 

\begin{lem} \label{lem_independent_perturbation}
Subject to Assumptions \ref{connected_root_node} and \ref{assum_trans_stability}, the multi-agent system \eqref{closed_loop} is stable independent of unitary phase perturbations $\mathrm{\bm{\Delta}}_p$ in the feedback path if and only if
\begin{equation}
\label{condition_independent_perturbation}
\bar{\sigma} (\mathbf{G}_p(j\omega_p)) <1, \quad \forall \omega_p > 0, \quad \forall p = 2, \dots, N
\end{equation}
where $\sigma_i (\mathbf{G}_p)$ are the singular values of the transfer function matrix $\mathbf{G}_p$, $\bar{\sigma}(\mathbf{G}_p) = \max \sigma_i (\mathbf{G}_p)$ and  $\underline{\sigma}(\mathbf{G}_p) = \min \sigma_i (\mathbf{G}_p)$ .
\end{lem}

\begin{IEEEproof}
Let us assume condition \eqref{condition_independent_perturbation} holds. Now, we can write
\begin{equation}
\bar{\sigma} \left((j \omega_p \mathbf{I} - \mathbf{A})^{-1} c \lambda_p \mathbf{B} \mathbf{K}  \right) < 1, \quad \forall \omega_p > 0.
\end{equation}

For unitary phase perturbation $e^{\bm{\Sigma}_p}$ in the feedback path, we have $\bar{\sigma} \left((j \omega_p \mathbf{I} - \mathbf{A})^{-1} c\mathbf{B} \mathbf{K} \lambda_p e^{\bm{\Sigma}_p} \right)<1$ which also can be expressed as
\begin{equation}
\label{sigma_bar_second_expression}
 \bar{\sigma} \left((j \omega_p \mathbf{I} - \mathbf{A})^{-1} c \lambda_p \mathbf{B} \mathbf{K} \bm{\Delta}_p \right) < 1
\end{equation}
where, $\bm{\Delta}_p=e^{\bm{\Sigma}_p}$ is unitary. It is straightforward to see that if condition \eqref{sigma_bar_second_expression} holds, then
\begin{equation*}
\text{det} \left(\mathbf{I} + (j \omega_p \mathbf{I} - \mathbf{A})^{-1} c \lambda_p \mathbf{B} \mathbf{K} \bm{\Delta}_p \right) \neq 0, \quad \forall \omega_p > 0
\end{equation*}
or equivalently, 
\begin{equation*}
\text{det} \left(\mathbf{I} + \mathbf{G}_p(j \omega_p) \bm{\Delta}_p \right) \neq 0, \quad \forall \omega_p > 0
\end{equation*}
i.e. the characteristic polynomial of the system \eqref{perturbed_trans_leaderless} does not intersect the imaginary axis and the system is stable independent of unitary phase perturbation. Moreover, from Lemma \ref{lem_alternate_stability} and Remark \ref{rem_alternate_stability}, the multi-agent system \eqref{closed_loop} is stable independent of unitary phase perturbations in the feedback path. The proof for the sufficiency part is completed. 

To establish the necessity, assume that  $\bar{\sigma} \left(\mathbf{G}_p(j\omega_{c_p}) \right) = \bar{\sigma} \left((j \omega_p \mathbf{I} - \mathbf{A})^{-1}  c \lambda_p  \mathbf{B} \mathbf{K} \right) = 1$, for some $\omega_{c_p} > 0$, for any $ p = 2, \dots, N$. This implies that there exists some unitary $\bm{\Delta}_p = e^{\bm{\Sigma}_p}$ such that $\text{det} \left(\mathbf{I} + (j \omega_{c_p} \mathbf{I} - \mathbf{A})^{-1}  c \lambda_p  \mathbf{B} \mathbf{K} \bm{\Delta}_p \right) = \text{det} \left(\mathbf{I} + \mathbf{G}_p(j \omega_p) \bm{\Delta}_p \right) = 0$ and from Lemma \ref{lem_instability_overall} the multi-agent system \eqref{closed_loop} becomes unstable. Let us now consider a case when, $\bar{\sigma} (\mathbf{G}_p(j\omega_p)) = \bar{\sigma} \left((j \omega_p \mathbf{I} - \mathbf{A})^{-1}  c \lambda_p \mathbf{B} \mathbf{K} \right) > 1$, for some $\omega_p > 0$. Since, $\bar{\sigma} \left(\mathbf{G}_p(j\omega_p)\right)$ is a continuous function of $\omega_p$, there exists some $\omega_{c_p} \in (\omega_p, \infty)$, such that $\bar{\sigma} \left(\mathbf{G}_p(j\omega_{c_p}) \right)= 1$ and the multi-agent system \eqref{closed_loop} is unstable. 
\end{IEEEproof}

\begin{rem}
\label{rem_sigma_bar}
Note that, if $\bar{\sigma} \left(\mathbf{G}_p(j\omega_{p}) \right) = 1$, there exists a unit vector $\mathbf{z}_p$ such that $\| \mathbf{G_p}(j \omega_p) \mathbf{z}_p\| = 1$. The proof of which is trivial and well known.
\end{rem}

Now let us define a set $\Omega_p = \{\omega_p |~ \underline{\sigma} (\mathbf{G}_p(j \omega_p)) \leq 1 \leq \bar{\sigma} (\mathbf{G}_p(j \omega_p)) \}$ for all $p = 2, \dots, N$. The cardinality of set $\Omega_p$ is denoted as $n_{\Omega_p}$. 

\subsubsection{Stability of multi-agent system dependent on unitary phase perturbations}
If the conditions highlighted by Lemma \ref{lem_independent_perturbation} are not satisfied, then there exists a unitary perturbation which destabilizes the multi-agent system. In this section, we provide the approach to find such perturbation and a computational framework to characterize the phase margin of the system.

\begin{lem}
\label{lem_unitary_mapping}
There exists a destabilizing unitary $\bm{\Delta}_p$ which is a mapping between two unit vectors, if and only if for any $p = 2, \dots, N$ the set $\Omega_p \neq \emptyset$.
\end{lem}

\begin{IEEEproof}
If the set $\Omega_p \neq \emptyset$ for any $p = 2, \dots, N$, there exists a $\omega_p$ and a unit vector $\mathbf{z}_p$ such that $\| \mathbf{G_p}(j \omega_p) \mathbf{z}_p \| = 1$ (see the proof of Lemma \ref{lem_independent_perturbation} and Remark \ref{rem_sigma_bar}). Let $\mathbf{v}_p = -\mathbf{G}_p (j \omega_p) \mathbf{z}_p$ and consider a unitary matrix $\mathrm{\bm{\Delta}}_p$ which maps $\mathbf{v}_p$ into $\mathbf{z}_p$ such that $\mathrm{\bm{\Delta}}_p \mathbf{v}_p = \mathbf{z}_p$. Since, $\mathbf{G}_p (j \omega_p) \mathrm{\bm{\Delta}}_p \mathbf{v}_p = \mathbf{G}_p (j \omega_p) \mathbf{z}_p = -\mathbf{v}_p$, one can write $(\mathbf{I}+\mathbf{G}_p (j \omega_p) \mathrm{\bm{\Delta}}_p)\mathbf{v}_p = 0$ which implies \text{det}$(\mathbf{I}+\mathbf{G}_p (j \omega_p) \mathrm{\bm{\Delta}}_p) = 0$ and the system $\mathbf{G}_p (j \omega_p)$ is unstable. This concludes the necessity.

Now, consider a destabilizing unitary $\mathrm{\bm{\Delta}}_p$  such that, $\text{det}(\mathbf{I} + \mathbf{G}_p (j \omega_p) \mathrm{\bm{\Delta}}_p) = 0$ and a unit vector $\mathbf{v}_p$ such that $(\mathbf{I}+\mathbf{G}_p (j \omega_p) \mathrm{\bm{\Delta}}_p)\mathbf{v}_p = 0$ and thus $\mathbf{G}_p (j \omega_p) \mathrm{\bm{\Delta}}_p \mathbf{v}_p = -\mathbf{v}_p$. Now let us assume $\mathrm{\bm{\Delta}}_p$ maps $\mathbf{v}_p$ into $\mathbf{z}_p$ such that $\mathbf{z}_p = \mathrm{\bm{\Delta}}_p \mathbf{v}_p$. As $\mathrm{\bm{\Delta}}_p$ is unitary and $\mathbf{v}_p$ is a unit vector, we have $\|\mathbf{z}_p\| = 1$.  So, we can write,
$\underline{\sigma}(\mathbf{G}_p) = \displaystyle \inf_{\|\mathbf{z}_p\| = 1} \|\mathbf{G}_p \mathbf{z}_p\| \leq 1$. Similarly, $1 \leq \displaystyle \sup_{\|\mathbf{z}_p\| = 1} \|\mathbf{G}_p \mathbf{z}_p\| = \bar{\sigma}(\mathbf{G}_p)$. Thus, the set $\Omega_p \neq \phi$. Hence, we have established sufficiency and necessity to the statement.
\end{IEEEproof}

\begin{thm}
\label{thm_phase_margin}
Suppose the Assumptions \ref{connected_root_node} and \ref{assum_trans_stability} hold. Let $\Pp$ be the set of all $p\subset \{2, \dots, N\}$ where $\Omega_p \neq \emptyset$. Then, the loop transfer function $\mathbf{G}_p(s)$ in \eqref{nominal_loop} is stable if the eigenvalues $\{\lambda_k(\bm{\Delta}_p)\}^n_{k=1}$ of unitary perturbation $\bm{\Delta}_p \in \C^n$ in the feedback path of $\mathbf{G}_p(s)$ for all $p \in \Pp$ satisfies $\max(|\text{Im}(\lambda_k(\bm{\Delta}_p))|) < \phi_p$ where
\begin{equation}
\phi_p = \min_{i=1, 2, \dots, n_{\Omega_p}} \{\phi_i \}
\label{pm_formula}
\end{equation}
and $\phi_i = \min \{\cos^{-1} \{\left\langle\mathbf{v}_p,\mathbf{z}_p\right\rangle\} \}$ with unit vectors $\mathbf{v}_p$ and $\mathbf{z}_p$ satisfying $\mathbf{v}_p = -\mathbf{G}_p(j \omega_p) \mathbf{z}_p$, for all $ \omega_p \in \Omega_p$. Moreover, the loop transfer function $\mathbf{G}_p(s)$ in \eqref{nominal_loop} is stable independent of unitary perturbation $\bm{\Delta}_p$ if $\Omega_p = \emptyset$ for all $p= 2, \dots, N$.
\end{thm}

\begin{IEEEproof}
From Lemma \ref{lem_unitary_mapping}, if for any $p = 2, \dots, N$ the set $\Omega_p \neq \emptyset$, then there exists an $\omega_p \in \Omega_p$ where the system destabilizes and a set of unit vectors $\mathbf{v}_p$ and $\mathbf{z}_p$ can be calculated that satisfies $\mathbf{v}_p = -\mathbf{G}_p(j \omega_p) \mathbf{z}_p$. Moreover, there also exists a destabilizing unitary perturbation (say $\bm{\Delta}_p^c$) that maps $\mathbf{v}_p$ to $\mathbf{z}_p$. 

For $\mathbf{G}_p(s)$ to be stable, phase of unitary $\mathrm{\bm{\Delta}}_p$ in the feedback path should be less than the smallest phase of destabilizing unitary perturbation $\bm{\Delta}_p^c$ that maps unit vector $\mathbf{v}_p$ to $\mathbf{z}_p$ for all $\omega_p \in \Omega_p$. Further, the angle between subspaces of $\C^n$ in which two unit vectors $\mathbf{v}_p$ and $\mathbf{z}_p$ lie is given by $\cos^{-1} \{\left\langle\mathbf{v}_p,\mathbf{z}_p\right\rangle\}$ \cite{galantai2006jordan}. Also, as $\bm{\Delta}_p^c$ is unitary, we can write $\left\langle\mathbf{v}_p,\mathbf{z}_p\right\rangle =  \left\langle\bm{\Delta}_p^c \mathbf{v}_p, \ \bm{\Delta}_p^c \mathbf{z}_p \right\rangle$. To that end, the phase of destabilizing $\bm{\Delta}_p^c$ which maps the two unitary vectors $\mathbf{v}_p$ and  $\mathbf{z}_p$  such that $\mathbf{z}_p = \bm{\Delta}_p^c \mathbf{v}_p$ is also $\cos^{-1} \{\left\langle\mathbf{v}_p,\mathbf{z}_p\right\rangle\}$. Henceforth, the smallest phase of destabilizing unitary perturbation for all $p \in \Pp$ can be obtained by minimizing $\cos^{-1} \{\left\langle\mathbf{v}_p,\mathbf{z}_p\right\rangle\}$ for all $\omega_p \in \Omega_p$ and  is given by
\begin{equation}
\phi_p = \min_{i=1, 2, \dots, n_{\Omega_p}} \{\phi_i \}, \quad \phi_i = \min \{\cos^{-1} \{\left\langle\mathbf{v}_p,\mathbf{z}_p\right\rangle\} \}.
\end{equation}

From Lemma \ref{lem_stab_boundary} and Remark \ref{rem_stab_boundary}, one can write phase of any unitary $\bm{\Delta}_p$ as $\max(|\text{Im}(\lambda_k(\bm{\Delta}_p))|)$. Therefore, for $\mathbf{G}_p(s)$ to be stable the eigenvalues of unitary perturbation $\bm{\Delta}_p$ in the feedback path of $\mathbf{G}_p(s)$ should satisfy $\max(|\text{Im}(\lambda_k(\bm{\Delta}_p))|) < \phi_p$. Further, if the set $\Omega_p = \emptyset$ for all $p=2, 3, \dots, N$, the multi-agent system remains stable independent of phase perturbation from Lemma \ref{lem_independent_perturbation}. This completes the proof.
\end{IEEEproof}

\begin{rem}
\label{rem_pm_formula}
Based on Theorem \ref{thm_phase_margin} and Remark \ref{rem_alternate_stability}, the phase margin of the multi-agent system can be calculated to be
\begin{equation}
\label{pm_formula_mas}
\phi^* =  \inf_{p \in \Pp} \left\lbrace \min_{\omega_p \in \Omega_p} \left\lbrace \min \{\cos^{-1} \{\left\langle\mathbf{v}_p,\mathbf{z}_p\right\rangle\} \right\rbrace \right\rbrace.
\end{equation}

Moreover, as cosine is a monotonically decreasing function in $[0, \pi]$, minimizing $\cos^{-1} \{\left\langle\mathbf{v}_p,\mathbf{z}_p\right\rangle\}$ is same as maximizing the inner product $\left\langle\mathbf{v}_p,\mathbf{z}_p\right\rangle$ satisfying $\mathbf{v}_p = -\mathbf{G}_p(j \omega_p) \mathbf{z}_p$, for all $ \omega_p \in \Omega_p$. 
\end{rem}

\subsubsection{Delay independent stability of multi-agent systems}

\begin{lem} \label{lem_independent_tdelay}
Subject to Assumptions \ref{connected_root_node} and \ref{assum_trans_stability}, the input delay multi-agent system \eqref{closed_loop_withdelay} is stable independent of delay if and only if
\begin{itemize}
\item[(i)] $\mathbf{A}$ is stable and
\item[(ii)] $\bar{\sigma} (\mathbf{G}_p(j\omega_p)) <1$, $\forall \omega_p > 0$, $\forall p = 2, \dots, N$.
\end{itemize}
\end{lem}

\begin{IEEEproof}
For the system to be stable independent of delay, it is necessary that it be stable for $\tau = \infty$, which requires condition (i) to hold (see \cite{gu2003stability}). Condition (ii) is neccessary and sufficient condition for the multi-agent system to be stable independent of unitary phase perturbations as discussed in Lemma \ref{lem_independent_perturbation}. As input delay links to a phase change with no gain change, condition (ii) is also necessary and sufficient for the system \eqref{closed_loop_withdelay} to be stable independent of delay.
\end{IEEEproof}

\subsubsection{Delay dependent stability of multi-agent systems}
The approach of characterizing the input delay margin of multi-agent delay system in this paper bears some similarity to that of ``frequency sweeping method" in the literature (see e.g., \cite{chen1995frequency, gu2003stability}). 
\begin{thm}
\label{prop_tdelay_margin}
Suppose the Assumptions \ref{connected_root_node} and \ref{assum_trans_stability} hold. Let $\Pp$ be the set of all $p\subset \{2, \dots, N\}$ where $\Omega_p \neq \emptyset$. Then, the input delay multi-agent system \eqref{closed_loop_withdelay} is stable for any $\tau \in [0, \tau^*)$ where
\begin{equation}
\tau^* = 
\begin{cases}
\displaystyle{\min_{p \in \Pp} \ \min_{1 \leq i \leq n_{\Omega_p}}} \ \dfrac{\phi_i}{\omega_i}, &\text{if } \ \Pp \neq \emptyset\\
\infty, \quad &\text{if } \ \Pp = \emptyset
\end{cases}
\label{td_formula}
\end{equation}
and  $\phi_i = \min \{\cos^{-1} \{\left\langle\mathbf{v}_p,\mathbf{z}_p\right\rangle\} \}$, $\omega_i = \text{argmin} \{\cos^{-1} \{\left\langle\mathbf{v}_p,\mathbf{z}_p\right\rangle\} \}$ with unit vectors $\mathbf{v}_p$ and $\mathbf{z}_p$ satisfying $\mathbf{v}_p = -\mathbf{G}_p(j \omega_p) \mathbf{z}_p$, for all $ \omega_p \in \Omega_p$.
\end{thm}

The proof follows from the proof of Theorem \ref{thm_phase_margin} and has been omitted for brevity. The sketch of the proof is as follows: since input delay can be linked to a unitary phase perturbation, once the phases $\phi_i$ are calculated, a set of delays can be calculated for each $\omega_p \in \Omega_p$ as $\tau_i = \dfrac{\phi_i}{\omega_i}$. Infimum of this set over all $\omega_p \in \Omega_p$ provides the upper limit of delay for the loop transfer function of $p^{\text{th}}$ input delayed system to remain stable, i.e. $\tau^*_p= \displaystyle\min_{1 \leq i \leq n_{\Omega_p}} \ \dfrac{\phi_i}{\omega_i}$ \cite{middleton2007achievable}. Moreover, from \ref{sec:mas_td}, one can establish $\tau^* = \displaystyle\min_{p \in \Pp} \tau^*_p$ such that the input delay multi-agent system \eqref{closed_loop_withdelay} is stable if $\tau \in [0, \tau^*)$. Further, if the set $\Omega_p = \emptyset$ for all $p=2, 3, \dots, N$, the system remains stable independent of delay from Lemma \ref{lem_independent_tdelay}.


\subsubsection{Computational Framework for Phase margin and Input delay margin} 
\label{sec:tdelay_calculation}
This section provides the computational framework to characterize the phase margin and input delay margin for multi-agent systems \eqref{closed_loop} and \eqref{closed_loop_withdelay}, respectively. In order to calculate the phase margin and input delay margin, one needs to find the set $\Omega_p$, for which it is necessary to find all $\omega_p>0$ such that $\bar{\sigma} (\mathbf{G}_p(j \omega_p)) \geq 1$ and $\underline{\sigma} (\mathbf{G}_p(j \omega_p)) \leq 1$, $\forall p =2, \dots, N$. The procedure to compute the set $\Omega_p$ is discussed in Procedure \ref{pm_td_algorithm}. As stated earlier, once the set $\Omega_p$ is calculated, the problem of calculating phase margin and input delay margin is equivalent to maximizing $<\mathbf{v}_p,\mathbf{z}_p>$ for all $\omega_p \in \Omega_p$ (see Remark \ref{rem_pm_formula}) which is same as maximizing $\left\langle\mathbf{v}_p,\mathbf{z}_p\right\rangle+\left\langle\mathbf{z}_p,\mathbf{v}_p\right\rangle= {\mathbf{v}^*_p\mathbf{z}_p+ \mathbf{z}^*_p \mathbf{v}_p}$. As $\mathbf{G}_p(j \omega_p) \mathbf{z}_p = -\mathbf{v}_p$, we can have
\begin{equation}
\begin{aligned}
\mathbf{v}^*_p\mathbf{z}_p+ \mathbf{z}^*_p \mathbf{v}_p =\ &  -\mathbf{z}^*_p \mathbf{G}_p(j \omega_p)^*\mathbf{z}_p - \mathbf{z}^*_p \mathbf{G}_p(j \omega_p) \mathbf{z}_p \\
=\ & -\mathbf{z}^*_p(\mathbf{G}_p(j \omega_p)^*+\mathbf{G}_p(j \omega_p))\mathbf{z}_p.  \label{vz}
\end{aligned}
\end{equation}

Now, maximizing $\mathbf{v}^*_p\mathbf{z}_p+ \mathbf{z}^*_p \mathbf{v}_p$ is equivalent to minimizing $\mathbf{z}^*_p(\mathbf{G}_p(j \omega_p)^*+\mathbf{G}_p(j \omega_p))\mathbf{z}_p$. Thus the problem of calculating phase margin is converted to a constrained minimization problem: minimize $\mathbf{z}^*_p(\mathbf{G}_p(j \omega_p)^*+\mathbf{G}_p(j \omega_p))\mathbf{z}_p$ such that $|\mathbf{v}_p| = |\mathbf{z}_p| =1$, $-\mathbf{G}_p(j \omega_p) \mathbf{z}_p = \mathbf{v}_p$ which can be further expressed as
\begin{equation}
\begin{aligned}
 &\text{minimize} \quad [\mathbf{z}_p^* (\mathbf{G}_p(j \omega_p) + \mathbf{G}_p(j \omega_p)^*) \mathbf{z}_p]\\
\text{subject to} \quad
 &\mathbf{z}^*_p\mathbf{z}_p = 1, \quad \mathbf{z}^*_p \ \mathbf{G}_p(j \omega_p)^* \ \mathbf{G}_p(j \omega_p) \ \mathbf{z}_p = 1.
\end{aligned} \label{opt_prob}
\end{equation}

Further discussion on the optimization problem is provided in Appendix \ref{appendix_II}. The complete procedure to compute the phase margin and input delay margin of the multi-agent system is discussed in Procedure \ref{pm_td_algorithm}.
\begin{algorithm}[t]
\caption{Computation of $\phi^*$ and $\tau^*$} 
\label{pm_td_algorithm}
\begin{algorithmic}[1]
\State  Calculation of set $\Omega_p$ for all $p =2, \dots, N$:
\begin{itemize}
\item[(i)] Solve det$\left(\mathbf{I}-\mathbf{G}_p(j \omega_p)^* \mathbf{G}_p(j \omega_p) \right) = 0$ for all real roots of $\omega_p$ and calculate the eigenvalues of $\mathbf{G}_p(j \omega_p)^* \mathbf{G}_p(j \omega_p)$ at each root $\omega_{p}$. Let $\omega_{k_p}$, $k \subset \{1,2, \dots \}$ denote all the real roots $\omega_{p}$.
\item[(ii)] Knowing the eigenvalues of $\mathbf{G}_p(j \omega_p)^* \mathbf{G}_p(j \omega_p)$ at each $\omega_k$ and at $0$ will enable one to determine if there exists a $\sigma(\mathbf{G}_p (j \omega_p)) \leq 1$ and a $\sigma(\mathbf{G}_p (j \omega_p)) \geq 1$ in the region $(\omega_{(k-1)_p}, \omega_{k_p}]$ with $\omega_{0_p} = 0$.
\item[(iii)] The set $\Omega_p$ can be obtained as $\Omega_p = \cup (\omega_{(k-1)_p}, \omega_{k_p})$.
\item[(iv)] If for any $p=2, 3, \dots, N$, $\sigma(\mathbf{G}_p (j \omega_p))$ does not span across 1, then $\Omega_p = \emptyset$.
\end{itemize} 
\State If for any $p=2, 3, \dots, N$, $\Omega_p \neq \emptyset$, solve optimization problem \eqref{opt_prob4} in Appendix  \ref{appendix_II} and compute $\mathbf{z}_p$ and $\mathbf{v}_p$ using \eqref{vz_after_opt} in Appendix  \ref{appendix_II}.
\State Compute $\phi^*$ and $\tau^*$ using \eqref{pm_formula_mas} and \eqref{td_formula}, respectively.
\end{algorithmic}
\end{algorithm}

\subsection{Gain Margin}
\label{sec:gm}
For the gain margin calculation, the gain information of $\bm{\Delta}_p$ is assumed to be contained in the positive definite Hermitian part $\mathbf{R}$ of the polar decomposition of $\bm{\Delta}_p$. The unitary part $\mathbf{U}$ is assumed to be lumped into the loop transfer function or assumed to be an identity matrix.

\begin{lem}
\label{lem_destab_gm_delta}
There exists a destabilizing positive definite Hermitian $\bm{\Delta}_p$ if and only if there exists an $\omega_p$ and a complex vector $\mathbf{z}_p$ such that
\begin{equation}
\left\langle \mathbf{G}_p(j \omega_p)\mathbf{z}_p, \mathbf{z}_p \right\rangle<0
\end{equation}
for any $p = 2, \dots, N$.
\end{lem}

\begin{IEEEproof}
Let $\mathbf{v}_p = -\mathbf{G}_p (j \omega_p) \mathbf{z}_p$. Now, if for any $p = 2, \dots, N$,  $\left\langle \mathbf{G}_p(j \omega_p)\mathbf{z}_p, \mathbf{z}_p \right\rangle<0$ implies $\mathbf{z}_p^* \mathbf{G}_p(j \omega_p)^* \mathbf{z}_p<0$, i.e.
\begin{equation}
\mathbf{v}_p^* \ \mathbf{z}_p >0 \label{vp_zp}
\end{equation}

Further if \eqref{vp_zp} holds, one can always find a positive definite Hermitian matrix $\mathrm{\bm{\Delta}}_p$ such that $\mathbf{z}_p = \mathrm{\bm{\Delta}}_p \mathbf{v}_p$, as discussed in Appendix \ref{appendix_I}.  Substituting $\mathbf{z}_p$ in \eqref{vp_zp}, we get
\begin{equation}
\mathbf{v}_p^* \ \mathrm{\bm{\Delta}}_p \mathbf{v}_p >0
\end{equation}

Moreover, since $\mathbf{G}_p (j \omega_p) \mathrm{\bm{\Delta}}_p \mathbf{v}_p = \mathbf{G}_p (j \omega_p) \mathbf{z}_p = -\mathbf{v}_p$, one can write $(\mathbf{I}+\mathbf{G}_p (j \omega_p) \mathrm{\bm{\Delta}}_p)\mathbf{v}_p = 0$ which implies \text{det}$(\mathbf{I}+\mathbf{G}_p (j \omega_p) \mathrm{\bm{\Delta}}_p) = 0$ and the system $\mathbf{G}_p (j \omega_p)$ is unstable. This concludes the necessity.

Now, consider a destabilizing positive definite Hermitian matrix $\mathrm{\bm{\Delta}}_p$  such that, $\text{det}(\mathbf{I} + \mathbf{G}_p (j \omega_p) \mathrm{\bm{\Delta}}_p) = 0$ and a unit vector $\mathbf{v}_p$ such that $(\mathbf{I}+\mathbf{G}_p (j \omega_p) \mathrm{\bm{\Delta}}_p)\mathbf{v}_p = 0$; thus, $\mathbf{G}_p (j \omega_p) \mathrm{\bm{\Delta}}_p \mathbf{v}_p = -\mathbf{v}_p$. Let us assume $\mathrm{\bm{\Delta}}_p$ maps $\mathbf{v}_p$ into $\mathbf{z}_p$ such that $\mathbf{z}_p = \mathrm{\bm{\Delta}}_p \mathbf{v}_p$, then $\mathbf{v}_p = -\mathbf{G}_p (j \omega_p) \mathbf{z}_p$. As $\mathrm{\bm{\Delta}}_p>0$, one can write
\begin{equation}
0 < \mathbf{v}^*_p \mathrm{\bm{\Delta}}_p \mathbf{v}_p  = \mathbf{v}^*_p \mathbf{z}_p = -\mathbf{z}_p^* \mathbf{G}_p(j \omega_p)^* \mathbf{z}_p 
\end{equation}
which leads to

\begin{equation}
\left\langle \mathbf{G}_p(j \omega_p)\mathbf{z}_p, \mathbf{z}_p \right\rangle < 0.
\end{equation}

 Hence, we have established sufficiency and necessity to the statement.
\end{IEEEproof}

\begin{rem}
\label{rem_real_pdh}
For positive definite Hermitian $\mathrm{\bm{\Delta}}_p$, $\mathbf{v}_p^* \ \mathrm{\bm{\Delta}}_p \mathbf{v}_p$ is always real and positive, i.e. $\mathbf{v}_p^* \mathbf{z}_p$ is also real and positive. Also, if $\mathbf{G}_p (j \omega_p) \mathbf{z}_p = -\mathbf{v}_p$, $\mathbf{z}_p^* \mathbf{G}_p(j \omega_p)^* \mathbf{z}_p$ is real, and thus $\mathbf{z}_p^* \mathbf{G}_p(j \omega_p)^* \mathbf{z}_p = \mathbf{z}_p^* \mathbf{G}_p(j \omega_p) \mathbf{z}_p$. Further, any positive definite Hermitian matrix $\bm{\Delta}_p$ can be written as $e^\mathbf{S}$. As $\lambda_k(\bm{\Delta}_p) = e^{\lambda_k(\mathbf{S})}$, for all $k = 1,2, \dots, n$, we define gain of $\bm{\Delta}_p$ as $\max|\lambda(\mathbf{S})|=\max|\ln(\lambda_k(\bm{\Delta}_p))|$. Note that, unlike in calculation of phase margin in section \ref{sec:pm_td}, $\mathbf{v}_{p}$ and $\mathbf{z}_{p}$ need not be unit vectors. 
\end{rem}

Now, let us define a set $\tilde{\Omega}_p = \{\omega_p |~ \left\langle \mathbf{G}_p(j \omega_p)\mathbf{z}_p, \mathbf{z}_p \right\rangle<0 \}$ for all $p = 2, \dots, N$. The cardinality of set $\tilde{\Omega}_p$ is denoted as $n_{\tilde{\Omega}_p}$.

\begin{cor}(Stability of multi-agent system independent of gain perturbations)
\label{cor_destab_gm}
If Assumptions \ref{connected_root_node} and \ref{assum_trans_stability} hold and the set $\tilde{\Omega}_p = \emptyset$ for all $p=2, 3, \dots, N$, the multi-agent system \eqref{closed_loop} remains stable independent of gain perturbation in the feedback path of each agents.
\end{cor}

\begin{thm}
\label{thm_gain_margin}
Suppose the Assumptions \ref{connected_root_node} and \ref{assum_trans_stability} hold. Let $\Pp$ be the set of all $p\subset \{2, \dots, N\}$ where $\tilde{\Omega}_p \neq \emptyset$. Then, the loop transfer function $\mathbf{G}_p(s)$ in \eqref{nominal_loop} is stable if any one of the following is satisfied:
\begin{itemize}
\item[(i)] Conditions of Corollary \ref{cor_destab_gm} hold, i.e., $\tilde{\Omega}_p = \emptyset$ for all $p = 2, \dots, N$.
\item[(ii)] if the eigenvalues $\{\lambda_k(\bm{\Delta}_p)\}^n_{k=1}$ of the positive definite Hermitian perturbation $\bm{\Delta}_p \in C^n$ in the feedback path of loop transfer function $\mathbf{G}_p(s)$ for all $p \in \Pp$ satisfy $\max|\ln(\lambda_k(\bm{\Delta}_p))| < g_p$ where
\begin{equation}
g_p \leq \displaystyle{\min_{1 \leq i \leq n_{\tilde{\Omega}_p}}} \ g_i \\
\end{equation}
and $g_i = \min \left\lbrace \cosh^{-1} \left[\dfrac{\mathbf{v}^*_p\mathbf{v}_p+\mathbf{z}^*_p \mathbf{z}_p}{2\mathbf{v}^*_p \mathbf{z}_p} \right] \right\rbrace$, with unit vectors $\mathbf{v}_p$ and $\mathbf{z}_p$ satisfying $\mathbf{v}_p = -\mathbf{G}_p(j \omega_p) \mathbf{z}_p$ for all $\omega_p \in \tilde{\Omega}_p$.
\end{itemize}
\end{thm}

\begin{IEEEproof}
Statement (i) follows from Corollary \ref{cor_destab_gm}. On the other hand, if for any $p = 2, \dots, N$ the set $\Omega_p \neq \emptyset$, then there exists a $\omega_p \in \tilde{\Omega}_p$ where the system destabilizes and a set of unit vectors $\mathbf{v}_p$ and $\mathbf{z}_p$ can be calculated that satisfies $\mathbf{v}_p = -\mathbf{G}_p(j \omega_p) \mathbf{z}_p$. Moreover, from Lemma \ref{lem_destab_gm_delta}, there also exists a destabilizing positive definite Hermitian perturbation (say $\bm{\Delta}_p^c$) that maps $\mathbf{v}_p$ to $\mathbf{z}_p$.

For $\mathbf{G}_p(s)$ to be stable, gain of positive definite Hermitian $\mathrm{\bm{\Delta}}_p$ in the feedback path should be less than the smallest gain of destabilizing positive definite Hermitian perturbation $\bm{\Delta}_p^c$ that maps unit vector $\mathbf{v}_p$ to $\mathbf{z}_p$ for all $\omega_p \in \tilde{\Omega}_p$. Further, from \cite{bar1991multivariable}, the gain between two complex vectors $\mathbf{z}_p$ and $\mathbf{v}_p$ is given by $ \left| \cosh^{-1} \left[\dfrac{\mathbf{v}^*_p\mathbf{v}_p+\mathbf{z}^*_p \mathbf{z}_p}{2\mathbf{v}^*_p \mathbf{z}_p} \right] \right|$. Moreover, gain between $\mathbf{z}_p$ and $\mathbf{v}_p$ is also the gain of the positive definite matrix that maps vectors $\mathbf{z}_p$ and $\mathbf{v}_p$. Note that $\mathbf{v}^*_p\mathbf{v}_p$ and $\mathbf{z}^*_p\mathbf{z}_p$ are positive and real, and from Remark \ref{rem_real_pdh},  $\mathbf{v}^*_p\mathbf{z}_p$ is also real and positive. Thus, $\cosh^{-1} \left[\dfrac{\mathbf{v}^*_p\mathbf{v}_p+\mathbf{z}^*_p \mathbf{z}_p}{2\mathbf{v}^*_p \mathbf{z}_p} \right]$ is real and positive.  Henceforth, the smallest gain of destabilizing positive definite Hermitian perturbation for all $p \in \Pp$ can be obtained by minimizing $  \cosh^{-1} \left[\dfrac{\mathbf{v}^*_p\mathbf{v}_p+\mathbf{z}^*_p \mathbf{z}_p}{2\mathbf{v}^*_p \mathbf{z}_p} \right]$ and  is given by
\begin{equation}
g_p = \min_{i=1, 2, \dots, n_{\Omega_p}} \{g_i \}, \quad g_i = \min \left\lbrace \cosh^{-1} \left[\dfrac{\mathbf{v}^*_p\mathbf{v}_p+\mathbf{z}^*_p \mathbf{z}_p}{2\mathbf{v}^*_p \mathbf{z}_p} \right] \right\rbrace.
\end{equation}

From Remark \ref{rem_real_pdh}, gain of $\bm{\Delta}_p$ is $\max|\ln(\lambda_k(\bm{\Delta}_p))|$. Therefore, for $\mathbf{G}_p(s)$ to be stable the eigenvalues of positive definite Hermitian perturbation $\bm{\Delta}_p$ in the feedback path of $\mathbf{G}_p(s)$ should satisfy $\max|\ln(\lambda_k(\bm{\Delta}_p))| < g_p$. This completes the proof.
\end{IEEEproof}

\begin{rem}
\label{rem_gain_margin}
Based on Theorem \ref{thm_gain_margin} and Remark \ref{rem_alternate_stability}, the gain margin of the multi-agent system which is the gain of the positive definite Hermitian matrix in the feedback path of each agents can be calculated to be
\begin{equation}
\label{gm_formula_mas}
g^* =  \displaystyle\begin{cases} \displaystyle\inf_{p \in \Pp} \left\lbrace \min_{\omega_p \in \Omega_p} \left\lbrace \cosh^{-1} \left[\dfrac{\mathbf{v}^*_p\mathbf{v}_p+\mathbf{z}^*_p \mathbf{z}_p}{2\mathbf{v}^*_p \mathbf{z}_p} \right] \right\rbrace \right\rbrace, \quad &\text{if } \Pp \neq \emptyset\\
\infty, \quad &\text{if } \Pp = \emptyset 
\end{cases}
\end{equation}

Further, if $\Omega_p = \emptyset$, it is straightforward to see that the multi-agent system is stable if the eigenvalues $\{\lambda_k\}^n_{k=1}$ of positive definite Hermitian matrix in the feedback path of all agents satisfy $\lambda_k \in \left[e^{-g^*}, \ e^{g^*} \right]$. Since, eigenvalues $\{\lambda_k\}^n_{k=1}$ and singular values $\{\sigma_k\}^n_{k=1}$ of a positive definite Hermitian matrices are equivalent, we have $\sigma_k \in \left[e^{-g^*}, \ e^{g^*} \right]$. 
\end{rem}

\subsubsection{Computational framework to calculate gain margin}
In order to calculate the gain margin, it is necessary to calculate the set $\tilde{\Omega}_p$. As stated in Lemma \ref{lem_destab_gm_delta} and Remark \ref{rem_real_pdh}, for a destabilizing positive definite Hermitian matrix in the feedback path of $\mathbf{G}_p(j \omega_p)$ to exist for any $p = 2, \dots, N$, $\mathbf{v}_p^* \ \mathbf{z}_p$ must be real and positive such that $\mathbf{v}_p = -\mathbf{G}_p (j \omega_p) \mathbf{z}_p$. This leads to
\begin{equation}
\begin{aligned}
\text{Re} (\mathbf{v}_p^* \ \mathbf{z}_p) =& \dfrac{1}{2} (\mathbf{v}_p^* \ \mathbf{z}_p + \mathbf{z}_p^* \ \mathbf{v}_p) \\
=\ &  \mathbf{z}_p^* \left[ -\dfrac{1}{2} \left(\mathbf{G}_p(j \omega_p)^*+\mathbf{G}_p(j \omega_p) \right)\right] \mathbf{z}_p\\
=\ &\mathbf{z}_p^*\mathbf{X}_p(j \omega_p) \mathbf{z}_p >0 \label{real_vz}
\end{aligned}
\end{equation}
and
\begin{equation}
\begin{aligned}
\text{Im} (\mathbf{v}_p^* \ \mathbf{z}_p) =\ & -\dfrac{1}{2} j(\mathbf{v}_p^* \ \mathbf{z}_p - \mathbf{z}_p^* \ \mathbf{v}_p)\\ =& \mathbf{z}_p^* \left[ -\dfrac{1}{2} j \left(\mathbf{G}_p(j \omega_p)-\mathbf{G}_p(j \omega_p)^*\right)\right] \mathbf{z}_p\\
=\ &\mathbf{z}_p^*\mathbf{Y}_p(j \omega_p) \mathbf{z}_p =0. \label{imag_vz}
\end{aligned}
\end{equation}

Note that, both $\mathbf{X}_p(j \omega_p)$ and $\mathbf{Y}_p(j \omega_p)$ in \eqref{real_vz} and \eqref{imag_vz} are Hermitian matrices which can be obtained by decomposing $\mathbf{G}_p(j \omega_p)$ as $\mathbf{G}_p(j \omega_p)$ = $\mathbf{X}+j\mathbf{Y}$ such that
\begin{equation}
\begin{aligned}
\mathbf{X}=\ & \dfrac{1}{2} \left(\mathbf{G}_p(j \omega_p)+\mathbf{G}_p(j \omega_p)^*\right), \quad \text{and}\\
\mathbf{Y}=\ & -\dfrac{1}{2} j\left(\mathbf{G}_p(j \omega_p)-\mathbf{G}_p(j \omega_p)^*\right)
\end{aligned}
\label{complex_decomposition}
\end{equation}

Now, for $\omega_p \in \tilde{\Omega}_p$, $\mathbf{X}_p(j \omega_p)$ needs to be positive definite and $\mathbf{Y}_p(j \omega_p)$ needs to be have an eigenvalue equal to zero simultaneously at $\omega_p$. The detailed procedure to calculate the set $\tilde{\Omega}_p$ is discussed in Procedure \ref{gm_algorithm}.

Once the set $\tilde{\Omega}_p$ is computed, we need to minimize $\dfrac{\mathbf{v}^*_p\mathbf{v}_p+\mathbf{z}^*_p \mathbf{z}_p}{\mathbf{v}^*_p \mathbf{z}_p}$ at each $\omega_p \in \tilde{\Omega}_{p}$. Since $\cosh$ is a monotonically increasing function on $[0, \infty)$, minimizing $\cosh^{-1} \left[\dfrac{\mathbf{v}^*_p\mathbf{v}_p+\mathbf{z}^*_p \mathbf{z}_p}{\mathbf{v}^*_p \mathbf{z}_p}\right]$ is same as minimizing $\dfrac{\mathbf{v}^*_p\mathbf{v}_p+\mathbf{z}^*_p \mathbf{z}_p}{\mathbf{v}^*_p \mathbf{z}_p}$.  Let us choose a normalization constant $\gamma^2 = \mathbf{v}^*_p\mathbf{z}_p$ such that $\tilde{\mathbf{v}}_p=\dfrac{1}{\gamma}\mathbf{v}_p$ and $\tilde{\mathbf{z}}_p=\dfrac{1}{\gamma}\mathbf{z}_p$. Note that $\tilde{\mathbf{v}}^*_p \tilde{\mathbf{z}}_p = 1$. With necessary simplifications, the minimization problem to calculate minimum gain destabilizing $\mathrm{\bm{\Delta}}_i$ can be written as
\begin{equation}
\begin{aligned}
 \text{minimize} \quad &\tilde{\mathbf{v}}^*_p \tilde{\mathbf{v}}_p + \tilde{\mathbf{z}}^*_p \tilde{\mathbf{z}}_p\\
\text{subject to} \quad
 &\tilde{\mathbf{v}}^*_p \tilde{\mathbf{z}}_p = 1\\
&\tilde{\mathbf{v}}_p = -\mathbf{G}_{p}(j \omega_i) \tilde{\mathbf{z}}_i
\end{aligned}
\end{equation}

 As $\mathbf{G}_p(j \omega_p) \mathbf{z}_p = -\mathbf{v}_p$ also implies $\mathbf{G}_p(j \omega_p) \tilde{\mathbf{z}}_p = -\tilde{\mathbf{v}}_p$, we have 
\begin{equation}
\begin{aligned}
\tilde{\mathbf{v}}^*_p \tilde{\mathbf{v}}_p+ \tilde{\mathbf{z}}^*_p \tilde{\mathbf{z}}_p =\ & \tilde{\mathbf{z}}^*_p \mathbf{G}_p(j \omega_p)^* \mathbf{G}_p(j \omega_p) \tilde{\mathbf{z}}_p + \tilde{\mathbf{z}}^*_p \mathbf{z}_p \\
 =\ & \tilde{\mathbf{z}}^*_p(\mathbf{G}_p(j \omega_p)^*\mathbf{G}_p(j \omega_p) + \mathbf{I}_n)\tilde{\mathbf{z}}_p
 \end{aligned}
\end{equation} 
and
\begin{equation}
 \tilde{\mathbf{v}}^*_p\tilde{\mathbf{z}}_p = \mathbf{z}^*_p \mathbf{G}_p(j \omega_p)^*  \tilde{\mathbf{z}}_p
 = \tilde{\mathbf{z}}^*_p \mathbf{G}_p(j \omega_p) \tilde{\mathbf{z}}_p
\end{equation} 

The last equality follows from Lemma \ref{rem_real_pdh}. Now the problem of calculating gain margin is converted to a constrained minimization problem: 
\begin{equation}
\begin{aligned}
\text{minimize} \quad &[\tilde{\mathbf{z}}^*_p(\mathbf{G}_p(j \omega_p)^*\mathbf{G}_p(j \omega_p) + \mathbf{I}_n)\tilde{\mathbf{z}}_p]\\
\text{subject to} \quad
&\text{Re} [\tilde{\mathbf{z}}^*_p \mathbf{G}_p(j \omega_p) \tilde{\mathbf{z}}_p] = 1\\
&\text{Im} [\tilde{\mathbf{z}}^*_p \mathbf{G}_p(j \omega_p) \tilde{\mathbf{z}}_p] = 0
\end{aligned} \label{opt_prob_gm}
\end{equation}

Further discussion on the optimization problem is provided in Appendix \ref{appendix_III}. The procedure to compute the gain margin of the multi-agent system is discussed in Procedure \ref{gm_algorithm}.

\begin{algorithm}[t]
\caption{Computation of $g^*$} 
\label{gm_algorithm}
\begin{algorithmic}[1]
\State  Calculation of set $\tilde{\Omega}_p$ for all $p =2, \dots, N$:
\begin{itemize}
\item[(i)] Find $\mathbf{X}(j \omega_p)$ and $\mathbf{Y}(j \omega_p)$ from \eqref{complex_decomposition} Solve det$\left(\mathbf{Y}(j \omega_p) \right) = 0$ for all real roots of $\omega_p$. Let $\omega_{k}$, $k \subset \{1,2, \dots \}$ denote all real roots $\omega_{p}$.
\item[(ii)] Calculate the eigenvalues of $\mathbf{Y}(j \omega_p)$ at each  $\omega_{k}$ and at 0. 
\item[(iii)] If for any $\omega \in (\omega_{(k-1)_p}, \omega_{k_p}]$  with $\omega_{0_p} = 0$, $\lambda_{\max} (\mathbf{Y}(j \omega)) \lambda_{\min} (\mathbf{Y}(j \omega)) \leq 0$ and $\mathbf{X}_p(j \omega)$ is positive semidefinite, then $(\omega_{(k-1)_p}, \omega_{k_p}] \subset \tilde{\Omega}_p$.
\item[(iii)] The set $\tilde{\Omega}_p$ can be obtained as $\tilde{\Omega}_p = \cup (\omega_{(k-1)_p}, \omega_k)$.
\item[(iv)] If for any $p=2, 3, \dots, N$, condition (iii) does not hold, then $\Omega_p = \emptyset$.
\end{itemize} 
\State If for any $p=2, 3, \dots, N$, $\Omega_p \neq \emptyset$, solve optimization problem \eqref{opt_prob4_gm} in Appendix  \ref{appendix_III} and compute $\tilde{\mathbf{z}}_p$ and $\tilde{\mathbf{\mathbf{v}}}_p$ using \eqref{vz_after_opt_gm} in Appendix  \ref{appendix_III}.
\State Compute $g^*$ using \eqref{gm_formula_mas}.
\end{algorithmic}
\end{algorithm}


\section{Simulation Results}
\label{sec:results}

To demonstrate the preceeding analysis, we consider a multi-agent system with following system matrices \cite{li2009consensus}:
\begin{equation}
\mathbf{A} = \begin{bmatrix}
-2 & 2\\
-1 & 1
\end{bmatrix}, \qquad \mathbf{B} = \begin{bmatrix}
1\\
0
\end{bmatrix}.
\label{system_dyn_simulation}
\end{equation}

The choice of $\mathbf{A}$ and $\mathbf{B}$ satisfies Assumption \ref{assum_stab}. Let us now choose a stabilizing feedback gain, $\mathbf{K} = \begin{bmatrix}
-2 & -0.5
\end{bmatrix}$ such that $\mathbf{A}-\mathbf{B} \mathbf{K}$ is Hurwitz. We consider a network of 3 agents with following graph Laplacian matrix,
\begin{equation}
\mathbf{L} = \begin{bmatrix}
0 & 0 & 0\\
-1 & 2 & -1\\
0 & -1 & 1
\end{bmatrix}.
\label{laplacian_simulation}
\end{equation} 

In order to calculate the value of coupling gain $c$ we follow the procedure described in \ref{sec:mas_gen} which is taken from \cite{li2009consensus, li2014distributed}. The characteristic polynomial of $\mathbf{A}- \sigma \mathbf{B} \mathbf{K}$ is calculated to be $p(s) = s^2+(1-2x-j2y)s+(5/2)x+j(5/2)y$ with $\sigma = x+jy$. From Lemma 4 of \cite{li2009consensus}, $\mathbf{A}- \sigma \mathbf{B} \mathbf{K}$ is stable if and only if $1-2x>0$ and $(25/2)(1-2x)^2 x^2-5y^2(1-2x)-(25/4)y^2 > 0$; which describes the consensus region $\mathcal{S}(x,y) = \{x+jy \mid x<0.5; (25/2)(1-2x)^2 x^2-5y^2(1-2x)-(25/4)y^2 > 0 \}$. From \cite{li2014distributed}, for the agents to reach consensus, the coupling gain $c$ is to be selected such that $c \lambda_p$, $p=2, 3, \dots, N$ belong to the consensus region $\mathcal{S}(x,y)$ where $\lambda_p$ are eigenvalues of Laplacian matrix. The non-zero eigenvalues of the Laplacian matrix are calculated to be: $\lambda_2 = 0.3820$ and $\lambda_3 = 2.6180$. Thus, $c<0.1910$ guarantees the consensus. For the simulation, we consider $c$ to be 0.15.

 Based on the framework provided in section \ref{sec:pm_td} and \ref{sec:gm}, the phase and gain margins are calculated to be $\phi^* = 0.1820$ radians and $g^* = 0.4025$,  respectively. In other words, any unitary matrix whose phase is less than $0.1820$ radians in the feedback path will not destabilize the system. From Lemma \ref{lem_stab_boundary} the stabilizing boundary of phase is symmetric about the origin; thus the overall phase margin of the multi-agent system is calculated to be $[-0.1820, 0.1820]$ radians. Moreover, as stated in Remark \ref{rem_gain_margin}, any positive definite Hermitian matrix in the feedback path of loop transfer function of each agents whose singular values lie within $\sigma^*=\left[e^{-g^*}, e^{g^*} \right] = [0.6686, 1.4956]$ would guarantee the stability of multi-agent system.
 
 Further, we compare our results with the conventional \emph{disk-based gain} and \emph{disk-based phase margins} that have been widely utilized in the literature as robustness measure of a general MIMO system and can be obtained from the sensitivity and complimentary sensitivity functions of the system \cite{schug:hal-01636859, furuta1987pole,blight1994practical}. To compare the conservativeness and accuracy, the obtained gain and phase margins from the proposed approach are compared with the disk-based gain and disk-based phase margins obtained from sensitivity and complimentary sensitivity functions of the multi-agent system. The disk-based gain margin in terms of singular values of perturbation matrix is calculated to be $\tilde{\sigma}^* = [0.5143, 1.0820]$ and the disk-based phase margin is calculated to be $[-0.0788,0.0788]$ radians.

To verify the accuracy of the proposed framework, we construct a matrix $\bm{\Delta} \in \C^2$ which can be polar decomposed as follows
\begin{equation}
\bm{\Delta} = \mathbf{R} \mathbf{U}
\end{equation}
where $\mathbf{R}$ is the positive definite Hermitian and $\mathbf{U}$ is a unitary matrix. As discussed in Lemma \ref{lem_stab_boundary}, we can construct $\mathbf{U}$ as
\begin{equation}
\mathbf{U} =  \mathbf{P} e^{ \{\text{diag}(j \phi_1, j\phi_2) \}} \mathbf{P}^*
\end{equation}
where $\mathbf{P}$ is any unitary matrix. To construct $\mathbf{U}$, we choose a unitary $\mathbf{P} = \begin{bmatrix}
\cos(0.2) & -\sin(0.2)\\
\sin(0.2) & \cos(0.2)
\end{bmatrix}$ and the phases of the unitary matrix to be $\phi_1 = 0.18$ radians and $\phi_2 = 0.16$ radians. Let $\mathbf{R} = \begin{bmatrix}
1 & -0.15\\
-0.15 & 1
\end{bmatrix}$ whose singular values are $ \sigma = [0.85,1.15]$. This yields
\begin{equation}
\bm{\Delta}=\begin{bmatrix}
0.9841+j0.1777 & -0.1487- j 0.0202\\
-0.1483- j 0.0229 & 0.9872+j0.1595
\end{bmatrix}. \label{delta_sim}
\end{equation}

Note that $\sigma^* \ni \sigma \notin \tilde{\sigma}^*$, and $\phi^* > \phi_1 > \tilde{\phi}^*$ and $\phi^* > \phi_2 > \tilde{\phi}^*$. Figure \ref{fig:state_trajectory_gm_pm} shows the states of agents with $\bm{\Delta}$ from \eqref{delta_sim} whose gain and phase are within the margins provided by the proposed approach but not within the margins provided by the disk-based margin. 

\begin{figure}[thpb]
\centering
      \includegraphics[width = \columnwidth]{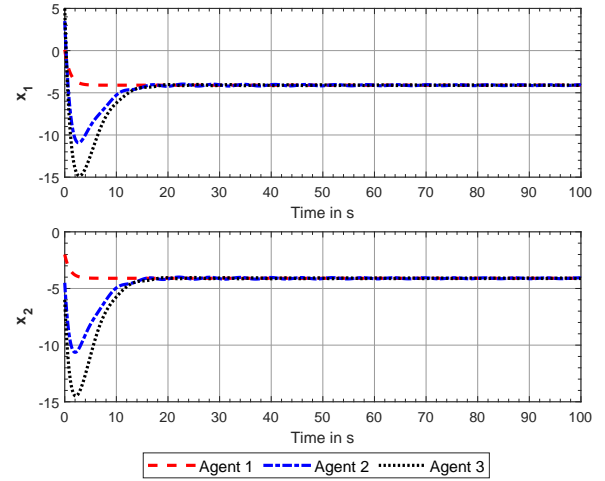}
      \caption{State trajectories of agents with feedback perturbation $\bm{\Delta}$ from \eqref{delta_sim}}
	  \label{fig:state_trajectory_gm_pm}
 \end{figure} 
 
Moreover, for the multi-agent system with input delay and with same system matrices as in \eqref{system_dyn_simulation} and graph Laplacian as in \eqref{laplacian_simulation}, the time delay margin is calculated to be $\tau^*=0.1978$ seconds. Figure \ref{fig:state_trajectory_gm_pm} shows the states of agents with a delay of $\tau = 0.18$ seconds in the inputs of three agents.

\begin{figure}[thpb]
\centering
      \includegraphics[width = \columnwidth]{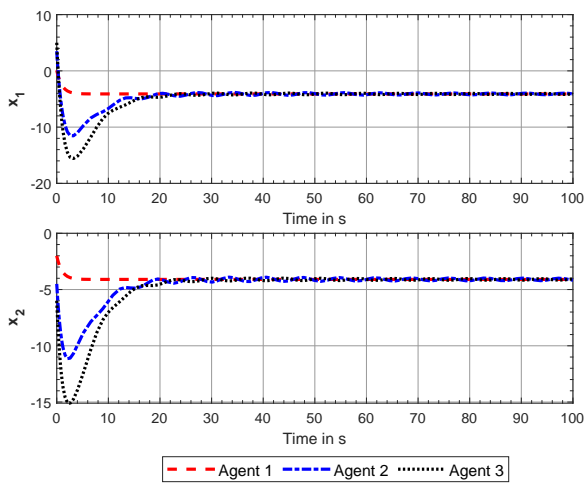}
      \caption{State trajectories of agents with $\tau = 0.18$ seconds}
	  \label{fig:state_trajectory_gm_pm}
 \end{figure} 
 
To illustrate the effectiveness of the proposed approach, we use different graph structures for the agents with the same system matrices as in \eqref{system_dyn_simulation} and with feedback gain matrix $\mathbf{K} = \begin{bmatrix}
-2 & -0.5
\end{bmatrix}$. For a directed cycle among 4 agents with $\mathbf{L} = \begin{bmatrix}
1 & 0 & 0 & -1\\
-1 & 1 & 0 & 0\\
0 & -1 & 1 & 0\\
0 & 0 & -1 & 1
\end{bmatrix}$, the distributed consensus protocol achieves consensus for any $c<0.5$. With $c = 0.15$, we compute gain margin and phase margin from the proposed approach to be $[0.3355, 2.9805]$ and $[-0.7995, 0.7995]$ radians, respectively. On the other hand, the disk-based gain and disk-based phase margins are computed to be $[0.676, 1.4792]$ and $[-0.3819, 0.3819]$ radians, respectively. Moreover, the input delay margin from the proposed approach is computed to be $2.05091$ seconds. 

Further, for an undirected cycle among 5 agents with $\mathbf{L} = \begin{bmatrix}
2&-1&0&0&-1\\
-1&2&-1&0&0\\
0&-1&2&-1&0\\
0&0&-1&2&-1\\
-1&0&0&-1&2\\
\end{bmatrix}$, the distributed consensus protocol achieves consensus for any $c<0.1382$. With $c = 0.12$, the gain margin and phase margin from the proposed approach is calculated to be $[0.6673, 1.4986]$ and $[-0.1066, 0.1066]$ radians, respectively. On the other hand, the disk-based gain and disk-based phase margins are computed to be $[0.6980, 1.0472]$ and $[-0.0461, 0.0461]$ radians, respectively. Moreover, the input delay margin from the proposed approach is computed to be $0.1066$ seconds. \emph{To that end, the proposed approach provides less conservative and accurate gain and phase margins within which the multi-agent system remains stable and achieves consensus, compared to disk-based gain and phase margins.}


\section{Conclusion}
\label{sec:conclusion}
In this paper, we have studied the stability of the multi-agent system under gain, phase, and input delay perturbations where each agent in the graph-based interconnection network is a linear time-invariant multi-input multi-output system. Based on the consensus protocol under a static graph communication topology, we provide a computational strategy to compute the gain, phase and input delay margins for multi-agent systems using the approach of multiplicative perturbation. Conditions for the gain, phase and delay independent stability of multi-agent system are discussed. To illustrate the effectiveness of the proposed framework, a numerical example with various graph structures was presented which depicted the lower conservativeness of the proposed approach as compared to disk-based gain and phase margins.


\bibliographystyle{IEEEtran}        
\bibliography{Bhusal_Subbarao_GainPhaseDelay_2020}

\begin{thebibliography}{10}
\providecommand{\url}[1]{#1}
\csname url@samestyle\endcsname
\providecommand{\newblock}{\relax}
\providecommand{\bibinfo}[2]{#2}
\providecommand{\BIBentrySTDinterwordspacing}{\spaceskip=0pt\relax}
\providecommand{\BIBentryALTinterwordstretchfactor}{4}
\providecommand{\BIBentryALTinterwordspacing}{\spaceskip=\fontdimen2\font plus
\BIBentryALTinterwordstretchfactor\fontdimen3\font minus
  \fontdimen4\font\relax}
\providecommand{\BIBforeignlanguage}[2]{{%
\expandafter\ifx\csname l@#1\endcsname\relax
\typeout{** WARNING: IEEEtran.bst: No hyphenation pattern has been}%
\typeout{** loaded for the language `#1'. Using the pattern for}%
\typeout{** the default language instead.}%
\else
\language=\csname l@#1\endcsname
\fi
#2}}
\providecommand{\BIBdecl}{\relax}
\BIBdecl

\bibitem{dong2016time}
X.~Dong and G.~Hu, ``Time-varying formation control for general linear
  multi-agent systems with switching directed topologies,'' \emph{Automatica},
  vol.~73, pp. 47--55, 2016.

\bibitem{jafari2019biologically}
M.~Jafari and H.~Xu, ``A biologically-inspired distributed fault tolerant
  flocking control for multi-agent system in presence of uncertain dynamics and
  unknown disturbance,'' \emph{Engineering applications of artificial
  intelligence}, vol.~79, pp. 1--12, 2019.

\bibitem{zhang2015cooperative}
Q.~Zhang, J.~Tao, F.~Yu, Y.~Li, H.~Sun, and W.~Xu, ``Cooperative solution of
  multi-{UAV} rendezvous problem with network restrictions,''
  \emph{Mathematical Problems in Engineering}, vol. 2015, 2015.

\bibitem{cai2014leader}
H.~Cai and J.~Huang, ``The leader-following attitude control of multiple rigid
  spacecraft systems,'' \emph{Automatica}, vol.~50, no.~4, pp. 1109--1115,
  2014.

\bibitem{safonov1981multiloop}
M.~Safonov and M.~Athans, ``A multiloop generalization of the circle criterion
  for stability margin analysis,'' \emph{IEEE Transactions on Automatic
  Control}, vol.~26, no.~2, pp. 415--422, 1981.

\bibitem{lehtomaki1981robustness}
N.~Lehtomaki, N.~Sandell, and M.~Athans, ``Robustness results in
  linear-quadratic gaussian based multivariable control designs,'' \emph{IEEE
  Transactions on Automatic Control}, vol.~26, no.~1, pp. 75--93, 1981.

\bibitem{nie2010exact}
Z.-Y. Nie, Q.-G. Wang, M.~Wu, and Y.~He, ``Exact computation of loop gain
  margins of multivariable feedback systems,'' \emph{Journal of Process
  Control}, vol.~20, no.~6, pp. 762--768, 2010.

\bibitem{wang2008loop}
Q.-G. Wang, Y.~He, Z.~Ye, C.~Lin, and C.~C. Hang, ``On loop phase margins of
  multivariable control systems,'' \emph{Journal of process Control}, vol.~18,
  no.~2, pp. 202--211, 2008.

\bibitem{tonetti2010limits}
S.~Tonetti and R.~M. Murray, ``Limits on the network sensitivity function for
  homogeneous multi-agent systems on a graph,'' in \emph{Proceedings of the
  2010 American Control Conference}.\hskip 1em plus 0.5em minus 0.4em\relax
  IEEE, 2010, pp. 3217--3222.

\bibitem{gattami2004frequency}
A.~Gattami and R.~Murray, ``A frequency domain condition for stability of
  interconnected mimo systems,'' in \emph{Proceedings of the 2004 American
  control conference}, vol.~4.\hskip 1em plus 0.5em minus 0.4em\relax IEEE,
  2004, pp. 3723--3728.

\bibitem{kim2017stability}
Y.~Kim, ``On the stability margin of networked dynamical systems,'' \emph{IEEE
  Transactions on Automatic Control}, vol.~62, no.~10, pp. 5451--5456, 2017.

\bibitem{cao2013overview}
Y.~Cao, W.~Yu, W.~Ren, and G.~Chen, ``An overview of recent progress in the
  study of distributed multi-agent coordination,'' \emph{IEEE Transactions on
  Industrial Informatics}, vol.~9, no.~1, pp. 427--438, 2013.

\bibitem{tian2008consensus}
Y.-P. Tian and C.-L. Liu, ``Consensus of multi-agent systems with diverse input
  and communication delays,'' \emph{IEEE Transactions on Automatic Control},
  vol.~53, no.~9, pp. 2122--2128, 2008.

\bibitem{xu2013input}
J.~Xu, H.~Zhang, and L.~Xie, ``Input delay margin for consensusability of
  multi-agent systems,'' \emph{Automatica}, vol.~49, no.~6, pp. 1816--1820,
  2013.

\bibitem{zhang2018synchronization}
M.~Zhang, A.~Saberi, and A.~A. Stoorvogel, ``Synchronization in a network of
  identical continuous-or discrete-time agents with unknown nonuniform constant
  input delay,'' \emph{International journal of robust and nonlinear control},
  vol.~28, no.~13, pp. 3959--3973, 2018.

\bibitem{olfati2004consensus}
R.~Olfati-Saber and R.~M. Murray, ``Consensus problems in networks of agents
  with switching topology and time-delays,'' \emph{IEEE Transactions on
  automatic control}, vol.~49, no.~9, pp. 1520--1533, 2004.

\bibitem{lin2008average}
P.~Lin and Y.~Jia, ``Average consensus in networks of multi-agents with both
  switching topology and coupling time-delay,'' \emph{Physica A: Statistical
  Mechanics and its Applications}, vol. 387, no.~1, pp. 303--313, 2008.

\bibitem{munz2010delay}
U.~M{\"u}nz, A.~Papachristodoulou, and F.~Allg{\"o}wer, ``Delay robustness in
  consensus problems,'' \emph{Automatica}, vol.~46, no.~8, pp. 1252--1265,
  2010.

\bibitem{zhang2019state}
M.~Zhang, A.~Saberi, A.~A. Stoorvogel, and Z.~Liu, ``State synchronization of a
  class of homogeneous linear multi-agent systems in the presence of unknown
  input delays via static protocols,'' \emph{European Journal of Control},
  vol.~47, pp. 20--29, 2019.

\bibitem{zhang2017distributed}
H.~Zhang, D.~Yue, W.~Zhao, S.~Hu, and C.~Dou, ``Distributed optimal consensus
  control for multiagent systems with input delay,'' \emph{IEEE transactions on
  cybernetics}, vol.~48, no.~6, pp. 1747--1759, 2017.

\bibitem{zhao2017guaranteed}
Y.~Zhao and W.~Zhang, ``Guaranteed cost consensus protocol design for linear
  multi-agent systems with sampled-data information: An input delay approach,''
  \emph{ISA transactions}, vol.~67, pp. 87--97, 2017.

\bibitem{bar1990phase}
J.~R. Bar-on and E.~A. Jonckheere, ``Phase margins for multivariable control
  systems,'' \emph{International Journal of control}, vol.~52, no.~2, pp.
  485--498, 1990.

\bibitem{bar1991multivariable}
------, ``Multivariable gain margin,'' \emph{International Journal of Control},
  vol.~54, no.~2, pp. 337--365, 1991.

\bibitem{yu2011second}
W.~Yu, W.~X. Zheng, G.~Chen, W.~Ren, and J.~Cao, ``Second-order consensus in
  multi-agent dynamical systems with sampled position data,''
  \emph{Automatica}, vol.~47, no.~7, pp. 1496--1503, 2011.

\bibitem{ren2007information}
W.~Ren, R.~W. Beard, and E.~M. Atkins, ``Information consensus in multivehicle
  cooperative control,'' \emph{IEEE Control systems magazine}, vol.~27, no.~2,
  pp. 71--82, 2007.

\bibitem{olfati2007consensus}
R.~Olfati-Saber, J.~A. Fax, and R.~M. Murray, ``Consensus and cooperation in
  networked multi-agent systems,'' \emph{Proceedings of the IEEE}, vol.~95,
  no.~1, pp. 215--233, 2007.

\bibitem{li2009consensus}
Z.~Li, Z.~Duan, G.~Chen, and L.~Huang, ``Consensus of multiagent systems and
  synchronization of complex networks: A unified viewpoint,'' \emph{IEEE
  Transactions on Circuits and Systems I: Regular Papers}, vol.~57, no.~1, pp.
  213--224, 2009.

\bibitem{li2014distributed}
Z.~Li and Z.~Duan, ``Distributed consensus protocol design for general linear
  multi-agent systems: a consensus region approach,'' \emph{IET Control Theory
  \& Applications}, vol.~8, no.~18, pp. 2145--2161, 2014.

\bibitem{zhang2011optimal}
H.~Zhang, F.~L. Lewis, and A.~Das, ``Optimal design for synchronization of
  cooperative systems: state feedback, observer and output feedback,''
  \emph{IEEE Transactions on Automatic Control}, vol.~56, no.~8, pp.
  1948--1952, 2011.

\bibitem{gantmakher1959theory}
F.~R. Gantmakher, \emph{The theory of matrices}.\hskip 1em plus 0.5em minus
  0.4em\relax AMS Chelsea Publishing, 1959, vol.~1.

\bibitem{zielinski1995polar}
P.~Zieli{\'n}ski and K.~Zi{{e}}tak, ``The polar decomposition-properties,
  applications and algorithms,'' \emph{Mathematica Applicanda}, vol.~24,
  no.~38, 1995.

\bibitem{chen2019phase}
W.~Chen, D.~Wang, S.~Z. Khong, and L.~Qiu, ``Phase analysis of {MIMO} {LTI}
  systems,'' in \emph{2019 IEEE 58th Conference on Decision and Control
  (CDC)}.\hskip 1em plus 0.5em minus 0.4em\relax IEEE, 2019, pp. 6062--6067.

\bibitem{galantai2006jordan}
A.~Gal{\'a}ntai and C.~J. Heged{\H{u}}s, ``Jordan's principal angles in complex
  vector spaces,'' \emph{Numerical Linear Algebra with Applications}, vol.~13,
  no.~7, pp. 589--598, 2006.

\bibitem{gu2003stability}
K.~Gu, J.~Chen, and V.~L. Kharitonov, \emph{Stability of time-delay
  systems}.\hskip 1em plus 0.5em minus 0.4em\relax Springer Science \& Business
  Media, 2003.

\bibitem{chen1995frequency}
J.~Chen and H.~A. Latchman, ``Frequency sweeping tests for stability
  independent of delay,'' \emph{IEEE Transactions on automatic control},
  vol.~40, no.~9, pp. 1640--1645, 1995.

\bibitem{middleton2007achievable}
R.~H. Middleton and D.~E. Miller, ``On the achievable delay margin using lti
  control for unstable plants,'' \emph{IEEE Transactions on Automatic Control},
  vol.~52, no.~7, pp. 1194--1207, 2007.

\bibitem{schug:hal-01636859}
A.~Schug, P.~Seiler, and H.~Pfifer, ``{Robustness Margins for Linear Parameter
  Varying Systems},'' \emph{{Aerospace Lab}}, no.~13, pp. pages 1--9, Nov.
  2017.

\bibitem{furuta1987pole}
K.~Furuta and S.~Kim, ``Pole assignment in a specified disk,'' \emph{IEEE
  Transactions on Automatic control}, vol.~32, no.~5, pp. 423--427, 1987.

\bibitem{blight1994practical}
J.~D. Blight, R.~Lane~Dailey, and D.~Gangsaas, ``Practical control law design
  for aircraft using multivariable techniques,'' \emph{International Journal of
  Control}, vol.~59, no.~1, pp. 93--137, 1994.

\bibitem{so2007approximating}
A.~M.-C. So, J.~Zhang, and Y.~Ye, ``On approximating complex quadratic
  optimization problems via semidefinite programming relaxations,''
  \emph{Mathematical Programming}, vol. 110, no.~1, pp. 93--110, 2007.

\bibitem{kuhn2014nonlinear}
H.~W. Kuhn and A.~W. Tucker, ``Nonlinear programming,'' in \emph{Traces and
  emergence of nonlinear programming}.\hskip 1em plus 0.5em minus 0.4em\relax
  Springer, 2014, pp. 247--258.

\bibitem{boyd2004convex}
S.~Boyd and L.~Vandenberghe, \emph{Convex optimization}.\hskip 1em plus 0.5em
  minus 0.4em\relax Cambridge university press, 2004.

\bibitem{guu1998quadratic}
S.~Guu and Y.~Liou, ``On a quadratic optimization problem with equality
  constraints,'' \emph{Journal of optimization theory and applications},
  vol.~98, no.~3, pp. 733--741, 1998.

\bibitem{bar1997global}
J.~Bar-On and K.~Grasse, ``Global optimization of a quadratic functional with
  quadratic equality constraints, part 2,'' \emph{Journal of Optimization
  Theory and Applications}, vol.~93, no.~3, pp. 547--556, 1997.

\end{thebibliography}
\appendices


\section{}
\label{appendix_I}
This section discusses a way to  finding a positive-definite Hermitian matrix $\mathbf{R}$ mapping $\mathbf{v} \in \C^n$ into $\mathbf{z} \in \C^n$, i.e.
\begin{equation}
\mathbf{z} = \mathbf{R} \mathbf{v} \label{eqn_zRv}
\end{equation}

Let the set of vectors $\{\mathbf{v}, \mathbf{z}, \mathbf{u}_1, \mathbf{u}_2, \dots, \mathbf{u}_{n-2} \}$ be a basis in $\C^n$. Given, a symmetric positive definite bilinear form $\langle \cdot, \cdot \rangle$  on finite-dimensional vector space, one can use the Gram-Schmidt orthogonalization process to find a orthonormal basis. Since $\langle \mathbf{v}, \mathbf{z}  \rangle = \mathbf{v}^{*} \mathbf{z} >0$, we can construct an orthonormal basis $\{\mathbf{q}_1, \mathbf{q}_2, \dots, \mathbf{q}_n \}$ as
\begin{equation}
\begin{aligned}
\mathbf{q}_1 =\ & \dfrac{\mathbf{v}}{\sqrt{\langle \mathbf{v}, \mathbf{v} \rangle}}; \quad \mathbf{q}_2 = \dfrac{\hat{\mathbf{q}}_2}{\sqrt{\langle \hat{\mathbf{q}}_2, \hat{\mathbf{q}}_2 \rangle}}, \quad \hat{\mathbf{q}}_2 =  \mathbf{z} - (\mathbf{q}^{*}_1 \mathbf{z}) \mathbf{q}_1;\\
\mathbf{q}_k =\ & \dfrac{\hat{\mathbf{q}}_k}{\sqrt{\langle \hat{\mathbf{q}}_k, \hat{\mathbf{q}}_k \rangle}}, \quad \hat{\mathbf{q}}_k =  \mathbf{u}_{k-2} - \sum^{k-1}_{i=1} (\mathbf{q}^{*}_i \mathbf{u}_{k-2}) \mathbf{q}_{i}, \quad k=3, \dots, n
\end{aligned}
\end{equation}

Here the matrix $\mathbf{Q} = [\mathbf{q}_1 \ \mathbf{q}_2\ \dots \ \mathbf{q}_n]$ is such that $\mathbf{Q} \mathbf{Q}^{*} = \mathbf{I}$. Pre-multiplying \eqref{eqn_zRv} with $\mathbf{Q}^*$, we obtain 
\begin{equation}
\mathbf{Q}^{*} \mathbf{z} = \mathbf{Q}^{*}  \mathbf{R} \mathbf{Q} \mathbf{Q}^{*}\mathbf{v} \label{eqn_transzRv1}
\end{equation}

Let $\mathbf{Q}^*$ maps $\mathbf{v}$ into $\mathbf{e}_1$ and $\mathbf{Q}^*$ maps $\mathbf{z}$ into a linear combination of $\mathbf{e}_1$ and $\mathbf{e}_2$, such that $\mathbf{v}=\mathbf{Q}(\gamma \mathbf{e}_1)$  and $\mathbf{z} = \mathbf{Q} (\alpha \mathbf{e}_1+ \beta \mathbf{e}_2 )$, where $\gamma =\sqrt{\langle \mathbf{v}, \mathbf{v} \rangle}$,  $\alpha = (\mathbf{q}^*_1 \mathbf{z})$ and $\beta =  \sqrt{\langle \hat{\mathbf{q}}_2, \hat{\mathbf{q}}_2 \rangle}$. Substituting for $\mathbf{z}$ and $\mathbf{v}$ in \eqref{eqn_transzRv1}, we obtain:  $\mathbf{Q}^{*} \mathbf{Q} (\alpha \mathbf{e}_1+ \beta \mathbf{e}_2 )  =  (\mathbf{Q}^{*}  \mathbf{R} \mathbf{Q}) \mathbf{Q}^{*} \mathbf{Q}(\gamma \mathbf{e}_1)$. With $\hat{\mathbf{v}} = \gamma \mathbf{e}_1$ and $\hat{\mathbf{z}} =\alpha \mathbf{e}_1+ \beta \mathbf{e}_2$ , one can write
\begin{eqnarray}
\hat{\mathbf{z}}& =& \mathbf{P} \hat{\mathbf{v}} 
\end{eqnarray}
where $\mathbf{P} =  \mathbf{Q}^{*}  \mathbf{R} \mathbf{Q}$. As $\mathbf{Q}$ is orthonormal matrix, eigenvalues of both matrices $\mathbf{R}$ and $\mathbf{P}$ are same. Since, only the upper $2 \times 2$ block of $\mathbf{Q}$ is needed to map $\hat{\mathbf{v}} $ into $\hat{\mathbf{z}}$, let us decompose $\mathbf{P}$ as $
\mathbf{P}=\begin{bmatrix}
\hat{\mathbf{P}} & \mathbf{0}\\
\mathbf{0} & \mathbf{I}
\end{bmatrix}
$,
where $\hat{\mathbf{P}} \in \C^{2 \times 2}$. Now the problem of finding a positive-definite Hermitian matrix $\mathbf{R}$ mapping $\mathbf{v}$ into $\mathbf{z}$ is reduced to the problem of finding $\hat{\mathbf{P}}$. Moreover, for $\mathbf{R}$ to be positive definite, $\hat{\mathbf{P}}$ must be positive definite and must be of the form
\begin{equation}
\hat{\mathbf{P}} = \begin{bmatrix}
\dfrac{\alpha}{\gamma} & \dfrac{\beta}{\gamma}\\[1em]
\dfrac{\beta}{\gamma} & \mathbf{p}_{22}
\end{bmatrix}
\end{equation}
where $\mathbf{p}_{22}$ should be such that
\begin{eqnarray}
\mathbf{p}_{22} > \dfrac{\beta^2}{\gamma \alpha}
\end{eqnarray}

Finally, $\mathbf{R}$ can be computed as $\mathbf{R}=\mathbf{Q} \mathbf{P} \mathbf{Q}^{*}$. It should be noted that, such a positive definite Hermitian matrix $\mathbf{R}$ is not unique. 


\section{}
\label{appendix_II}
With $\mathbf{U} = \mathbf{G}_p(j \omega_p) + \mathbf{G}_p(j \omega_p)^*$, $\mathbf{V} = \mathbf{G}_p(j \omega_p)^* \ \mathbf{G}_p(j \omega_p)$ and $\mathbf{w} =\mathbf{z}_p$,  optimization problem in \eqref{opt_prob} can be rewritten as

\begin{equation}
\begin{aligned}
&\text{minimize} \quad [\mathbf{w}^* \mathbf{U} \mathbf{w}]\\
\text{subject to} \quad
&\mathbf{w}^*\mathbf{w} = 1\\
& \mathbf{w}^* \mathbf{V} \mathbf{w} = 1
\end{aligned} \label{opt_prob2}
\end{equation}

It is straightforward to show that the complex optimization problem in \eqref{opt_prob2} is equivalent to the following optimization problem from the work carried out in \cite{so2007approximating}:

\begin{equation}
\begin{aligned}
&\text{minimize} \quad  \left[(\mathbf{a}^{\text{T}},\mathbf{b}^{\text{T}})  \begin{pmatrix}
\text{Re}(\mathbf{U}) & \text{Im}(\mathbf{U})\\
-\text{Im}(\mathbf{U}) & \text{Re}(\mathbf{U})
\end{pmatrix} \begin{pmatrix}
\mathbf{a}\\
\mathbf{b}
\end{pmatrix}  \right]\\
\text{subject to} \quad
& \mathbf{a}_i^2+\mathbf{b}_i^2 = 1 \quad i=1,2,\dots, n\\
& (\mathbf{a}^{\text{T}},\mathbf{b}^{\text{T}})  \begin{pmatrix}
\text{Re}(\mathbf{V}) & \text{Im}(\mathbf{V})\\
-\text{Im}(\mathbf{V}) & \text{Re}(\mathbf{V}) \end{pmatrix} \begin{pmatrix}
\mathbf{a}\\
\mathbf{b}
\end{pmatrix}  = 1\\
& \mathbf{a}, \mathbf{b} \in \R^n
\end{aligned} \label{opt_prob3}
\end{equation}

Let $\mathbf{Q} = 
\begin{pmatrix}
\text{Re}(\mathbf{U}) & \text{Im}(\mathbf{U})\\
-\text{Im}(\mathbf{U}) & \text{Re}(\mathbf{U}) 
\end{pmatrix}$, $\mathbf{R} =
\begin{pmatrix}
\text{Re}(\mathbf{V}) & \text{Im}(\mathbf{V})\\
-\text{Im}(\mathbf{V}) & \text{Re}(\mathbf{V}) 
\end{pmatrix}$, $\mathbf{y}= \begin{pmatrix}
\mathbf{a}\\
\mathbf{b}
\end{pmatrix}$ in \eqref{opt_prob3}. Now, the transformed optimization problem becomes 

\begin{equation}
\begin{aligned}
&\text{minimize} \quad   \mathbf{y}^\text{T}  \mathbf{Q} \mathbf{y}\\
\text{subject to} \quad
& \mathbf{y}^\text{T} \mathbf{y}  = 1 \\
& \mathbf{y}^\text{T} \mathbf{R}  \mathbf{y}  = 1\\
& \mathbf{y} \in \R^{2n}
\end{aligned} \label{opt_prob4}
\end{equation}

The constrained optimization problem in \eqref{opt_prob4} is a set of quadratic optimization problems with nonlinear equality constraints which can be solved by solving the Karush-Kuhn-Tucker (KKT) optimality conditions \cite{kuhn2014nonlinear}. Moreover, the optimization problem in \eqref{opt_prob4} can be equivalently written as unconstrained minimization problem by defining the \textit{Lagrangian} as

\begin{equation}
\mathcal{L} (\mathbf{y}, \mu_1, \mu_2) = \mathbf{y}^{\text{T}} \mathbf{Q} \mathbf{y} + \mu_1 \left(\mathbf{y}^{\text{T}} \mathbf{y} - 1\right) +  \mu_2 \left(\mathbf{y}^{\text{T}} \mathbf{R} \mathbf{y} - 1\right) \label{Lag_prob}
\end{equation}
where, $\mu_1$ and $\mu_2$ are the scalar Lagrange multipliers associated with the equality constraints \cite{boyd2004convex}. Let ($\mathbf{y}^o$, $\mu^o_1$, $\mu^o_2$)  be the optimal solution to the optimization problem. Since $\mathbf{y}^*$ minimizes $\mathcal{L} (\mathbf{y}, \mu^o_1, \mu^o_2)$ over $y$, its gradient must vanish at $\mathbf{y}^o$. Hence, the KKT conditions which are necessary for the optimality can be written as follows:
\begin{equation}
\begin{aligned}
\mathbf{y}^{o^\text{T}} \mathbf{y^o} - 1 =\ & 0\\
\mathbf{y}^{o^\text{T}} \mathbf{R} \mathbf{y^o} - 1 =\ &0\\
\mathbf{Q} \mathbf{y}^o + \mu^o_1 \mathbf{y}^o + \mu^o_2 \mathbf{R} \mathbf{y}^o =\ & 0 
\end{aligned} \label{KKT}
\end{equation}

The KKT optimality conditions are a set of $2n+2$ equations with $2n+2$ unknown variables. The optimal solution obtained from solving \eqref{KKT} system of equations is the global minima to the original problem in \eqref{opt_prob4} if following KKT sufficient optimality condition holds:
\begin{equation} 
\mathbf{Q} + \mu^o_1 \mathbf{I}_{2n} + \mu^o_2 \mathbf{R} \geq 0 \label{KKT_2}
\end{equation}

Thus, any numerical routine that can generate the local optimum ($\mathbf{y}^o$, $\mu^o_1$, $\mu^o_2$) by solving \eqref{KKT} and eventually satisfies \eqref{KKT_2} gives the global optimum $\mathbf{y}$. Further discussion on global optimization of the similar problem (quadratic objective function with quadratic equality constraints) can be found in \cite{guu1998quadratic,bar1997global}. Once vector $\mathbf{y}$ is obtained by solving the optimization problem \eqref{opt_prob4}, the vectors $\mathbf{a}\in \R^n$ and $\mathbf{b}\in \R^n$ can be calculated and, the vector $\mathbf{w}\in \C^n$ or equivalently $\mathbf{z}_p \in \C^n$ and $\mathbf{v}_p \in \C^n$ can be obtained as
\begin{equation}
\begin{aligned}
\mathbf{z}_p =\ & \mathbf{a} + j \mathbf{b} \\
\mathbf{v}_p =\ & -\mathbf{G}_p(j \omega_p) \mathbf{z}_p 
\end{aligned} 
\label{vz_after_opt}
\end{equation} 


\section{}
\label{appendix_III}
With $\mathbf{U} = \mathbf{G}_p(j \omega_p)$, $\mathbf{V} = \mathbf{G}_p(j \omega_p)^* \ \mathbf{G}_p(j \omega_p)$ and $\mathbf{w} =\tilde{\mathbf{z}}_p$,  optimization problem in \eqref{opt_prob_gm} can be rewritten as

\begin{equation}
\begin{aligned}
&\text{minimize} \quad  [\mathbf{w}^* \mathbf{V} \mathbf{w}]\\
\text{subject to} \quad
&\mathbf{w}^* \mathbf{U} \mathbf{w} = 1\\
& \mathbf{w}^* (j\mathbf{I}) \mathbf{U} \mathbf{w} = 0
\end{aligned} \label{opt_prob2_gm}
\end{equation}

It is straightforward to show that the complex optimization problem in \eqref{opt_prob2_gm} is equivalent to the following optimization problem.

\begin{equation}
\begin{aligned}
&\text{minimize} \quad   \left[(\mathbf{a}^{\text{T}},\mathbf{b}^{\text{T}})  \left[\begin{pmatrix}
\text{Re}(\mathbf{V}) & \text{Im}(\mathbf{V})\\
-\text{Im}(\mathbf{V}) & \text{Re}(\mathbf{V})
\end{pmatrix} + \mathbf{I}_{2n} \right] \begin{pmatrix}
\mathbf{a}\\
\mathbf{b}
\end{pmatrix}  \right]\\
&\text{subject to} \\
& (\mathbf{a}^{\text{T}},\mathbf{b}^{\text{T}})  \begin{pmatrix}
\text{Re}(\mathbf{U}) & \text{Im}(\mathbf{U})\\
-\text{Im}(\mathbf{U}) & \text{Re}(\mathbf{U}) \end{pmatrix} \begin{pmatrix}
\mathbf{a}\\
\mathbf{b}
\end{pmatrix}  = 1\\
& (\mathbf{a}^{\text{T}},\mathbf{b}^{\text{T}}) \begin{pmatrix}
\text{Re}(j\mathbf{I}) & \text{Im}(j\mathbf{I})\\
-\text{Im}(j\mathbf{I}) & \text{Re}(j\mathbf{I}) \end{pmatrix} \begin{pmatrix}
\text{Re}(\mathbf{U}) & \text{Im}(\mathbf{U})\\
-\text{Im}(\mathbf{U}) & \text{Re}(\mathbf{U}) \end{pmatrix} \begin{pmatrix}
\mathbf{a}\\
\mathbf{b}
\end{pmatrix}  = 0\\
& \mathbf{a}, \mathbf{b} \in \R^n
\end{aligned} \label{opt_prob3_gm}
\end{equation}

Now, to further simplify \eqref{opt_prob3_gm}, let $\mathbf{Q} = 
\begin{pmatrix}
\text{Re}(\mathbf{U}) & \text{Im}(\mathbf{U})\\
-\text{Im}(\mathbf{U}) & \text{Re}(\mathbf{U}) 
\end{pmatrix}$, $\mathbf{R} =
\begin{pmatrix}
\text{Re}(\mathbf{V}) & \text{Im}(\mathbf{V})\\
-\text{Im}(\mathbf{V}) & \text{Re}(\mathbf{V}) 
\end{pmatrix}$, $\mathbf{y}= \begin{pmatrix}
\mathbf{a}\\
\mathbf{b}
\end{pmatrix}$, $\mathbf{J}=\begin{pmatrix}
\text{Re}(j\mathbf{I}) & \text{Im}(j\mathbf{I})\\
-\text{Im}(j\mathbf{I}) & \text{Re}(j\mathbf{I}) \end{pmatrix}$. To that end, optimization problem in \eqref{opt_prob3_gm} becomes 

\begin{equation}
\begin{aligned}
&\text{minimize} \quad   \mathbf{y}^\text{T}  \left(\mathbf{R} + \mathbf{I}_{2n} \right) \mathbf{y}\\
\text{subject to} \quad
& \mathbf{y}^\text{T}  \mathbf{Q} \mathbf{y}  = 1 \\
& \mathbf{y}^\text{T} \mathbf{J} \mathbf{Q}  \mathbf{y}  = 0\\
& \mathbf{y} \in \R^{2n}
\end{aligned} \label{opt_prob4_gm}
\end{equation}

The optimization problem in \eqref{opt_prob4_gm} can be equivalently written as unconstrained minimization problem by defining the \textit{Lagrangian} as in \eqref{Lag_prob} and similar KKT conditions as in \eqref{KKT} and \eqref{KKT_2} can be derived. Once vector $\mathbf{y}$ is obtained by solving the optimization problem \eqref{opt_prob4_gm}, the vectors $\mathbf{a}\in \R^n$ and $\mathbf{b}\in \R^n$ can be calculated and, the vector $\mathbf{w}\in \C^n$ or equivalently $\tilde{\mathbf{z}}_p \in \C^n$ and $\tilde{\mathbf{v}}_p \in \C^n$ can be obtained as
\begin{equation}
\begin{aligned}
\tilde{\mathbf{z}}_p =\ & \mathbf{a} + j \mathbf{b} \\
\tilde{\mathbf{v}}_p =\ & -\mathbf{G}_p(j \omega_p) \tilde{\mathbf{z}}_p
\end{aligned} 
\label{vz_after_opt_gm}
\end{equation} 

\end{document}